\documentclass[a4paper,11pt,twoside]{amsart}

\title{The Universal Functorial Equivariant Lefschetz Invariant}
\author{Julia Weber}


\address{Julia Weber\newline Max-Planck-Institut f\"ur Mathematik\newline Vivatsgasse 7\newline D-53111 Bonn\newline Germany\newline email: jweber@mpim-bonn.mpg.de}

\usepackage{bbm}
\usepackage{amsmath}          
\usepackage{amssymb} \usepackage{amsthm} \usepackage{xspace}
\usepackage[all]{xy}\CompileMatrices
\usepackage{epsfig}           
\usepackage{amsfonts} \usepackage{amscd} \usepackage{color}
\usepackage{stmaryrd}

\SelectTips{cm}{} \makeindex

 \DeclareMathOperator{\fp}{fp}
 \DeclareMathOperator{\End}{End}

\DeclareMathOperator{\Mor}{Mor} \DeclareMathOperator{\MMod}{Mod}
\DeclareMathOperator{\Aut}{Aut} \DeclareMathOperator{\Ob}{Ob}
\DeclareMathOperator{\cst}{cst} \DeclareMathOperator{\pt}{pt}
\DeclareMathOperator{\Deg}{Deg} \DeclareMathOperator{\ch}{ch}
\DeclareMathOperator{\Fix}{Fix} \DeclareMathOperator{\inc}{inc}
\DeclareMathOperator{\tr}{tr} \DeclareMathOperator{\iind}{ind}
\DeclareMathOperator{\rres}{res} \DeclareMathOperator{\id}{id}
\DeclareMathOperator{\Id}{Id} 
\DeclareMathOperator{\ff}{f\/f} \DeclareMathOperator{\el}{el}
\DeclareMathOperator{\cone}{cone} \DeclareMathOperator{\coker}{coker}

\DeclareMathOperator{\consub}{consub}
 
 \DeclareMathOperator{\Is}{Is}
\DeclareMathOperator{\rk}{rk} \DeclareMathOperator{\pr}{pr}
 \DeclareMathOperator{\incl}{incl}
 
 \DeclareMathOperator{\Split}{Split}

\DeclareMathOperator{\colim}{colim}

\newcommand{\C}{\ensuremath{\mathcal{C}}\xspace}
\newcommand{\Z}{\ensuremath{\mathbb{Z}}\xspace}
\newcommand{\Q}{\ensuremath{\mathbb{Q}}\xspace}

\newcommand{\N}{\ensuremath{\mathbb{N}}\xspace}

\newcommand{\oD}{\ensuremath{\text{\emph{\r{D}}}}\xspace}
\newcommand{\ind}[1]{\iind_{#1}}
\newcommand{\rrres}[1]{\rres_{#1}}
\newcommand{\encat}[2]{{#1}\text{\normalfont -end}_{\text{\normalfont
f\/f}R{#2}}} \newcommand{\encatch}[2]{{#1}\text{\normalfont
-end}_{\text{\normalfont f\/f}R{#2}\text{\normalfont -ch}}}
 \newcommand{\encatp}[2]{{#1}\text{\normalfont
-end}_{\text{\normalfont fp}R{#2}}}
\newcommand{\encatpz}[2]{{#1}\text{\normalfont
-end}_{\text{\normalfont fp}\Z{#2}}}
   \newcommand{\bet}{\emph}
\newcommand{\Ab}{\mathcal{A}b} \newcommand{\CW}{\mathcal{C}W}
\newcommand{\Ch}{\mathcal{C}h} 
 
\newcommand{\RMod}{R \text{-$\MMod$}} \newcommand{\RCh}{R
\text{-$\Ch$}}

\newcommand{\EndoGCWcat}{\End(G\text{-}\CW_{\fp})}
\newcommand{\GCWcat}{G\text{-}\CW_{\fp}}

\theoremstyle{definition} \newtheorem{Def}{Definition}[section]
\newtheorem{Ex}[Def]{Example} 
\newtheorem*{Def*}{Definition} \newtheorem{Rem}[Def]{Remark}

\theoremstyle{plain}
\newtheorem{Thm}[Def]{Theorem} \newtheorem{Cor}[Def]{Corollary}
\newtheorem{Prop}[Def]{Proposition} \newtheorem{Lem}[Def]{Lemma}
\newtheorem*{Thm*}{Theorem}
\newtheorem*{Thma}{Theorem~\ref{universalitytheorem}}
\newtheorem*{Thmb}{Theorem~\ref{0.2*}}

\begin{document}

\maketitle

\begin{abstract}
We introduce the universal functorial equivariant Lefschetz invariant
for endomorphisms of finite proper $G$-CW-complexes, where $G$ is a
discrete group. We use $K_0$ of the category of
``$\phi$-endomorphisms of finitely generated free
$R\Pi(G,X)$-modules''. We derive results about fixed points of
equivariant endomorphisms of cocompact proper smooth $G$-manifolds.
\end{abstract}

\section*{Introduction} 

The Lefschetz number is a classical invariant of algebraic
topology. Numerous generalizations have been studied, for example the
generalized Lefschetz invariant~\cite{reidemeister,wecken2} and the
Lefschetz zeta function~\cite{felshtyn-hill94,geoghegan-nicas94}.

An invariant
which maps to all of these generalizations and still has the
characteristic properties of the Lefschetz number, namely
homotopy invariance and additivity, has been constructed by
L\"uck~\cite{lueck99}. He defines a universal
functorial Lefschetz invariant $(U,u)$ for endomorphisms of finite
CW-complexes. Here, $U$ is a functor which
assigns an abelian group $U(X,f)$ to each endomorphism $f\colon X\to X$. The
function $u$ assigns an element $u(X,f)\in U(X,f)$ to the endomorphism~$f$.

The aim of this article is to generalize this construction to the
equivariant setting. On the one hand, this is interesting because it gives
finer invariants. If there is a $G$-action on a space, this extra
information is taken into account in the construction of the
invariant. On the other hand, this generalization enlarges the scope
of the invariant. If we have an infinite discrete group $G$ acting
properly on a finite $G$-CW-complex~$X$, the space $X$ seen as a
CW-complex is infinite. This situation cannot be treated by the
classical theory, but the equivariant version can be applied.

For compact Lie groups $G$, there are
equivariant versions of the generalized Lefschetz invariant
and related constructions~\cite{marzantowicz-prieto,wong91_equivnielsen}. For discrete groups $G$, equivariant Lefschetz
numbers have been defined~\cite{lueck-laitinen,lueck-rosenberg}, but
none of the generalizations have been studied. We
define an equivariant version of the universal invariant, which
entails equivariant versions of all intermediate generalizations of the
Lefschetz number.

We deal with discrete groups $G$ and finite
proper $G$-CW-complexes~$X$. To an equivariant endomorphism $f\colon
X\to X$, we associate an abelian group $U_G(X,f)$ and an element
$u_G(X,f)\in U_G(X,f)$. The pair $(U_G,u_G)$ is a functorial equivariant Lefschetz
invariant, i.e., it is additive, $G$-homotopy invariant
and compatible with the induction structure in $G$. We show that it
has a universal initial property:

\begin{Thm} \label{universalitytheorem}
The pair $(U_G,u_G)$ is the universal functorial equivariant
Lefschetz invariant for equivariant endomorphisms of finite proper
$G$-CW-complexes, for
discrete groups $G$.
\end{Thm}

The invariant $(U_G,u_G)$ contains much information, and the
 abelian group $U_G(X,f)$ can be quite large. We show that
 $U_G(X,f)$ has a direct sum decomposition
into summands corresponding to the fixed point sets $X^H$, for
subgroups $H\leq G$. This not only gives structural insight, but is also helpful for concrete calculations. The information contained in $u_G(f)$ splits up into the information given by the restrictions of $f$ to the pairs $(X^H,X^{>H})$, where $X^{>H}$ is the subset of $X^H$ of points with larger isotropy group than $H$.

The decomposition is obtained from a $K$-theoretic splitting theorem, Theorem~\ref{theorem431}. This splitting theorem is valid for all $K$-groups, not only for $K_0$. A more general version for endomorphisms categories of modules over categories can be formulated and applied to the study of $K$-theory of endomorphism categories of modules over group rings.

The study of Lefschetz invariants was motivated by interest in fixed points. As an application of our general constructions, we therefore
derive an invariant from $(U_G,u_G)$ extracting information about fixed points, the generalized equivariant Lefschetz invariant
$(\Lambda_G, \lambda_G)$. A trace map $\tr_G$ which maps $(U_G,u_G)$
to $(\Lambda_G, \lambda_G)$ is constructed. The use of generalized traces in equivariant fixed point theory already appears in~\cite{marzantowicz}, in the context of equivariant $K$-theory.

We prove the refined
equivariant Lefschetz fixed point theorem, which shows that the
generalized equivariant Lefschetz invariant $\lambda_G(f)$ is equal to
the sum of all fixed point contributions. Here, the fixed point data is collected in the generalized
local equivariant Lefschetz class $\lambda_G^{loc}(f)$.

\begin{Thm}\label{0.2*}
Let $G$ be a discrete group, let $M$ be a
  cocompact proper smooth $G$-manifold and let $f\colon M\to M$ be a
  $G$-equivariant endomorphism such that $\Fix(f)\cap \partial M =
  \emptyset$ and such that for every $x\in \Fix(f)$ the determinant of the
  map $(\id_{T_x M} - T_x f)$ is different from zero. Then
\[
\lambda_G(f)=\lambda^{loc}_G(f).
\]
\end{Thm}

There are further applications to the study of fixed points of equivariant endomorphisms of cocompact proper smooth $G$-manifolds.

Based on the generalized equivariant Lefschetz invariant $\lambda_G(f)$, we can introduce equivariant Nielsen invariants. These give lower bounds for the number of fixed point orbits in the $G$-homotopy class of $f$. Under mild hypotheses, these bounds are sharp. These results generalize results of Wong~\cite{wong93} to all discrete groups $G$.

In the same vein, we can define $G$-Jiang spaces and prove a converse to the equivariant Lefschetz fixed point theorem. In particular, under mild hypotheses, a $G$-equivariant endomorphism $f$ is $G$-homotopic to a fixed point free $G$-map if the generalized equivariant Lefschetz invariant~$\lambda_G(f)$ is zero.

These applications go beyond the scope of this present article, but they illustrate the interest of the invariants constructed here.

This paper is organized as follows:

1. The Algebraic Approach

2. The Geometric Approach

3. Universality

4. Splitting Results

5. The Generalized Equivariant Lefschetz Invariant

6. The Refined Equivariant Lefschetz Fixed Point Theorem

\section{The Algebraic Approach} \label{sec1} 

We start with purely algebraic considerations. Given a commutative
ring $R$, we define an abelian group $U(R,\Gamma,\phi)$ for any
EI-category $\Gamma$ with endofunctor $\phi$. We proceed to show
functoriality of this group in $(\Gamma,\phi)$ and analyze its
additive properties.

For us, the most important example of an EI-category is the
fundamental category $\Pi(G,X)$ of a topological space $X$ with an
action of a discrete group $G$~\cite[Definition~8.15]{lueck89}.

\begin{Def}
An \emph{EI-category} is a small category $\Gamma$ such that any endomorphism
in $\Gamma$ is an isomorphism.
An \emph{$R\Gamma$-module} is a contravariant functor
$M\colon \Gamma\to \RMod$. Given an endofunctor $\phi\colon \Gamma\to \Gamma$, we define a \emph{$\phi$-endomorphism} of an
$R\Gamma$-module $M$ to be a natural transformation $g\colon M\to M\circ
\phi$.
\end{Def}

We assemble all $\phi$-endomorphisms of finite free $R\Gamma$-modules in the category $\encat{\phi}{\Gamma}$. Morphisms from
$g\colon M\to M\circ \phi$ to $h\colon N\to N\circ \phi$ are natural
transformations $\tau\colon M\to N$ such that $h \tau=(\phi^*\tau)
g$. This is an exact category whose isomorphism classes of objects
form a set.

\begin{Def}
We define $U(R,\Gamma,\phi):= K_0(\encat{\phi}{\Gamma})$.
\end{Def}

Let $H\colon \Gamma_1 \to \Gamma_2$ be a covariant functor of
EI-categories with endofunctors $\phi$ and $\psi$ such that $\psi H=
H\phi$. We extend the covariant functor ``induction with $H$''~\cite[Definition~9.15]{lueck89},
defined on $R\Gamma_1$-modules, to
$\phi$-endomorphisms.

Induction $\ind{H}$ and restriction $\rrres{H}$ are adjoint
functors \cite[9.22]{lueck89}, so there are natural transformations
$\eta\colon \Id \to \rrres{H}\ind{H}$ and $\varepsilon\colon
\ind{H}\rrres{H} \to \Id$. We
define $H_*\colon \encat{\phi}{\Gamma_1}\to \encat{\psi}{\Gamma_2}$ on objects of $\encat{\phi}{\Gamma_1}$ to
be the composition $\varepsilon_{\rrres{\psi}\ind{H}}\circ
\ind{H}\rrres{\phi}(\eta) \circ \ind{H}$. By the
triangular identities~\cite[Theorem~IV.1.1]{maclane} we have $H_*(g)=\ind{H}(g)$ if $\phi=\Id$ and
$\psi=\Id$. We define $H_*$ on morphisms~$\tau$ by setting
$H_*(\tau):= \ind{H}(\tau)$. The functor $H_*$ is exact and thus
induces a group homomorphism $U(R,H):=K_0(H_*)\colon
U(R,\Gamma_1,\phi)\to U(R,\Gamma_2, \psi)$.

If we have another covariant functor $K\colon \Gamma_2\to \Gamma_3$ to
a category $\Gamma_3$ with endofunctor $\vartheta$ such that
$\vartheta K=K\psi$, we observe that $H_*G_*=(HG)_*$. We know even
more.

\begin{Prop}\label{equivofcat}
If $H\colon \Gamma_1\xrightarrow{\sim} \Gamma_2$ is an equivalence of
categories with endofunctors $\phi$ and $\psi$ such that $\psi
H=H\phi$, then the natural transformation $H_*\colon
\encat{\phi}{\Gamma_1} \to \encat{\psi}{\Gamma_2}$ is also an
equivalence of categories.
\end{Prop}

\begin{proof}
We use $\rrres{H}$ to define a functor $H^*\colon
\encat{\psi}{\Gamma_2}\to \encat{\phi}{\Gamma_1}$ by setting
$H^*(g):=\rrres{H}g$ on objects and $H^*(\tau):=\rrres{H}\tau$ on
morphisms. This is well-defined since $H$ is an equivalence of
categories. We now use the triangular identities to show that the
adjunction morphisms between $\ind{H}$ and $\rrres{H}$ induce natural
equivalences $\eta\colon \Id \to H^*H_*$ and $\varepsilon\colon
H_*H^*\to \Id$.
\end{proof}

We obtain a universal additive invariant for the objects of
$\encat{\phi}{\Gamma}$ by taking their representatives in
$U(R,\Gamma,\phi)$.

\begin{Def}\label{defaddinv}
An \emph{additive invariant} for the category $\encatch{\phi}{\Gamma}$ of $\phi$-endomorphisms
of finite free $R\Gamma$-chain complexes is a
pair $(A,a)$, where $A$ is an abelian group and $a\colon
\Ob(\encatch{\phi}{\Gamma})\to A, (g\colon C\to C\circ \phi)\mapsto
a(g)$ is a function, satisfying the following properties:
\begin{enumerate}
\item Additivity\\ For a short exact sequence $0\to f\to g\to h\to 0$,
given by a commutative diagram with exact rows
\[
\xymatrix{ 0 \ar[r] & C \ar[r] \ar[d]^{f} & D \ar[r] \ar[d]^{g} & E
\ar[r] \ar[d]^{h} & 0 \\ 0 \ar[r] & C \circ \phi \ar[r] & D \circ \phi
\ar[r] & E \circ \phi \ar[r] & 0 }
\]
where $C$, $D$ and $E$ are finite free $R\Gamma$-chain complexes, we
have
\[
a(f)-a(g) + a(h)=0.
\]
\item Homotopy invariance\\ Let $f, g\colon C\to C\circ \phi$ be
$R\Gamma$-chain maps of finite free $R\Gamma$-chain complexes. If $f$
and $g$ are $R\Gamma$-chain-homotopic (i.e., there is an $h\colon C\to
C\circ \phi$ of degree $+1$ such that $d h + h d = f-g$),
then $a(f)=a(g)$.
\item Elementary chain complexes\\ For a finite free $R\Gamma$-module
$F$ and $n\geq 1$ we have
\[
a\bigl(0\colon \el(F,n)\to \el(F,n)\circ\phi\bigr)=0,
\]
where $\el(F,n)$ denotes the elementary chain complex
concentrated in dimensions $n$ and $n-1$ and having $n$-th differential
$\Id\colon F\to F$.
\end{enumerate}
\end{Def}

\begin{Def}\label{defunivaddinv}
An additive invariant $(U,u)$ for
$\encatch{\phi}{\Gamma}$ is called \emph{universal} if and only if for
all additive invariants $(A,a)$ for the category
$\encatch{\phi}{\Gamma}$ there is exactly one group homomorphism
$\xi\colon U\to A$ such that $\xi(u(g))=a(g)$ for all $g\in
\encatch{\phi}{\Gamma}$.
\end{Def}
A finite free chain complex of $\phi$-endomorphisms maps into
$U(R,\Gamma,\phi)$ by $u\colon \{ \cdots \to g_n \to g_{n-1} \to
\cdots \} \mapsto \sum_{n\in \Z} (-1)^{n} [ g_n ] \in
U(R,\Gamma,\phi)$.

\begin{Thm}
The pair $\bigl(U(R,\Gamma,\phi),u\bigr)$ is the universal additive
invariant for $\encatch{\phi}{\Gamma}$.
\end{Thm}

\begin{proof}
Properties 1 and 3 are clear. In the proof of property 2, we proceed
analogously to L\"uck \cite[Theorem~1.4]{lueck99}, using the mapping
cone and the suspension of $C$. If $h\colon C\to C\circ \phi$ is an
$R\Gamma$-chain homotopy from $f$ to $g$, we obtain an exact sequence
$0\to f \to k \to \Sigma g \to 0$, where $k\colon \cone(C)\to
\cone(C)$ is given by $k_n:= \left( \begin{tiny} \begin{array}{cc}
g_{n-1} & 0 \\ h_{n-1} & f_n \end{array} \end{tiny} \right)$. Additivity yields $u(f)-u(g)=u(f)+u(\Sigma g)=u(k)$. The
$R\Gamma$-chain complex $\cone(C)$ is contractible. We
finally show that for contractible $C$ and for any $f\colon
C\to C\circ \phi$ we have $u(f\colon C\to C\circ \phi)=0\in
U(R,\Gamma,\phi)$.
\end{proof}

From the above proof, we can also conclude the following lemma.

\begin{Lem}\label{homequnull}
Let $f\colon C\to C\circ \phi$ and $g\colon D\to D\circ \phi$ be
$\phi$-endomorphisms of finite free $R\Gamma$-chain complexes. If
$h\colon C\to D$ is a chain homotopy equivalence such that $gh=(\phi^*h)
f$, then $u(f)=u(g)$ in $U(R,\Gamma,\phi)$.
\end{Lem}

\begin{proof}
We have a short exact sequence $0\to g \to k \to \Sigma f \to 0$,
where $ k_n=\left( \begin{tiny} \begin{array}{cc} -f_{n-1} & 0 \\
h_{n-1} & g_n \end{array} \end{tiny} \right).  $ We know that
$\cone(h)$ is contractible if and only if $h$ is a chain homotopy
equivalence. So $u(k)=0\in U(R,\Gamma,\phi)$,
whence $u(g)=u(f)$.
\end{proof}

If $R$ is a commutative ring which is not the zero ring, a finite free
$R$-module has a well-defined rank. One sees that the invariant
$\bigl(U(R,\Gamma,\phi),u\bigr)$ incorporates the Euler
characteristic. For a finite free $R\Gamma$-chain complex $C$, the
Euler characteristic is defined to be
$\xi(C):=\rk_{R\Gamma}([D])=\sum_{n\in \Z}(-1)^n\rk([C_n])\in
U(\Gamma):=\bigoplus_{\Is \Gamma}\Z$ \cite[Chapter~10]{lueck89}. There
is a natural transformation $s\colon U(\Gamma)\to U(R,\Gamma,\phi)$
mapping an element of $U(\Gamma)$ to the class of the $0$-map on the
corresponding free module. We have
$s(\chi(C))=u(0\colon C\to C\circ \phi)$.

\begin{Lem}
Let $v\colon C\to D$ and $f\colon C\to C\circ \phi$ be $R\Gamma$-chain
maps, where $C$ and $D$ are finite free $R\Gamma$-chain
complexes. Then
\[
u(f\circ v)+ s\bigl(\chi(D)\bigr)=u(\phi^* v \circ f ) +
s\bigl(\chi(C)\bigr).
\]
\end{Lem}

\begin{proof}
We have $\left( \begin{tiny} \begin{array}{cc} 1& \phi^* v \\ 0 & 1
\end{array} \end{tiny} \right) \left( \begin{tiny} \begin{array}{cc} 0
& 0 \\ f & f\circ v \end{array} \end{tiny} \right) = \left(
\begin{tiny} \begin{array}{cc} \phi^*v\circ f & 0 \\ f & 0 \end{array}
\end{tiny} \right) \left(
\begin{tiny} \begin{array}{cc}
1& v \\ 0 & 1 \end{array} \end{tiny} \right)$. Since the matrices with
$1$ on the diagonal are isomorphisms, we obtain from additivity $ u\Bigl(
\bigl( \begin{tiny} \begin{array}{cc} 0 & 0 \\ f & f\circ v
\end{array} \end{tiny} \bigr)\Bigr)= u\Bigl(\bigl( \begin{tiny}
\begin{array}{cc} \phi^*v\circ f & 0 \\ f & 0 \end{array} \end{tiny}
\bigr)\Bigr).  $ Again from additivity we derive
\begin{eqnarray*}
u(f\circ v) + u(0\colon D\to \phi^*D) & = & u\Bigl(\bigl( \begin{tiny}
\begin{array}{cc} 0 & 0 \\ f & f\circ v \end{array} \end{tiny}
\bigr)\Bigr)
\\
u(\phi^*v \circ f) + u(0\colon C\to \phi^*C) & = & u\Bigl(\bigl(
\begin{tiny} \begin{array}{cc} \phi^*v\circ f & 0 \\ f & 0 \end{array}
\end{tiny} \bigr)\Bigr).
\end{eqnarray*}
\vskip-4.5ex
\end{proof}

\section{The Geometric Approach} \label{secgeom} 

Let $G$ be a discrete group. Let $\GCWcat$ be the category of finite
proper $G$-CW-complexes. Let $\EndoGCWcat$ be the
category of $G$-equivariant endomorphisms of finite proper $G$-CW-complexes.

Given a commutative ring $R$, a finite proper $G$-CW-complex $X$ and a $G$-equivariant
cellular endomorphism $f\colon X\to X$, we want to introduce an
abelian group $U^R_G(X,f)$ and an invariant $u^R_G(X,f)\in
U^R_G(X,f)$. (The assumption that $f\colon X\to X$ is cellular can
finally be dropped because of homotopy invariance.)

In the non-equivariant case~\cite{lueck99}, the construction uses the
universal covering space $\widetilde{X}$ of the CW-complex
$X$. Choosing a basepoint $x$, the fundamental group $\pi_1(X,x)$ acts
on $\widetilde{X}$. One can lift $f\colon X\to X$ to a
$\phi$-equivariant map $\widetilde{f}\colon \widetilde{X}\to
\widetilde{X}$, where $\phi:=c_w\circ \pi_1(f,x) \colon \pi_1(X,x)\to
\pi_1(X,x)$ for a chosen path $w$ from $f(x)$ to $x$. We generalize
this non-equivariant construction to $G$-CW-complexes.

The natural generalization of the fundamental group is the fundamental
category $\Pi(G,X)$~\cite[Definition~8.15]{lueck89}. Objects in the fundamental category $\Pi(G,X)$
are $G$-maps $x\colon G/H\to X$ for $H\leq G$, morphisms are pairs
$(\sigma,[w])\colon x(H)\to y(K)$ consisting of a $G$-map
$\sigma\colon G/H\to G/K$ and a homotopy class $[w]$ relative
$G/H\times \partial I$ of $G$-maps $w\colon G/H\times I \to X$ with
$w_1=x$ and $w_0=y\circ \sigma$.  This
construction is functorial in $X$ (by composition of maps), in
particular the $G$-equivariant endomorphism $f\colon X\to X$ induces
an endofunctor $\phi:=\Pi(G,f)\colon \Pi(G,X)\to \Pi(G,X)$.

There is the contravariant universal covering functor
$\widetilde{X}\colon \Pi(G,X)\to \CW$ at hand \cite[Definition~8.22,
Proposition~8.33]{lueck89}, generalizing the universal
covering space. It maps $x\in \Pi(G,X)$ to $\widetilde{X^H(x)}$ and
$(\sigma,[w])\in \Mor(x(H),y(K))$ to the map $\widetilde{X^K(y)} \to
\widetilde{X^H(x)}$ induced by composition of morphisms. Here $X^H(x)$
denotes the connected component of the fixed point set $X^H$
containing the point $x(1H)$, for an object $x\colon G/H\to X$ of
$\Pi(G,X)$. Composing the universal covering functor $\widetilde{X}$
with the cellular chain complex functor $C^c\colon \CW \to \RCh$ one
obtains the cellular $R\Pi(G,X)$-chain complex $C^c(\widetilde{X})$ as a
contravariant functor $ C^c\circ
\widetilde{X}\colon \Pi(G,X)\xrightarrow{\widetilde{X}} \CW
\xrightarrow{C^c} \RCh$~\cite[Definition~8.37]{lueck89}. The functor $C^c\circ \widetilde{X}$ is a
finite free $R\Pi(G,X)$-chain complex~\cite[9.18]{lueck89}. The map $f\colon X\to X$ induces a $\phi$-endomorphism
$C^c(\widetilde{f})$ of the cellular $R\Pi(G,X)$-chain complex
$C^c\circ \widetilde{X}$.

\begin{Def}
We set $ U^R_G(X,f):=U\bigl(R,\Pi(G,X),\phi\bigr)$. We define the
element $u^R_G(X,f)\in U^R_G(X,f)$ to be $
u\bigl(C^c(\widetilde{f})\bigr)\in
U\bigl(R,\Pi(G,X),\phi\bigr)=U^R_G(X,f).  $
\end{Def}

Let $(X,f)$ and $(Y,g)$ be finite proper $G$-CW-complexes with
$G$-equivariant cellular endomorphisms $f$ and $g$. Let $h\colon X\to
Y$ be a $G$-equivariant cellular map such that $gh=hf$. Composition
with $h$ induces a functor $H:=\Pi(G,h)\colon \Pi(G,X)\to
\Pi(G,Y)$ between the fundamental categories. Setting $\phi:=\Pi(G,f)$
and $\psi:=\Pi(G,g)$, we have $\psi H=H\phi$. We obtain $U^R_G(h):=K_0(H_*)\colon U^R_G(X,f)\to U^R_G(Y,g)$, a group homomorphism. So
$U^R_G$ is a functor from $\EndoGCWcat$ to $\Ab$.

\begin{Lem}\label{cd}
If $h\colon (X,f)\to (Y,g)$ is a $G$-equivariant cellular map between
finite proper $G$-CW-complexes with endomorphisms such that $gh=hf$,
then $h$ induces a map $C^c(\widetilde{h})$ from
$\Pi(G,h)_*C^c(\widetilde{f})$ to $C^c(\widetilde{g})$.
\end{Lem}

\begin{proof}[Idea of proof]
We use adjointness of induction and restriction and the
triangular identities.
\end{proof}

Let $X$ be an $H$-CW-complex with $H$-equivariant endomorphism
$f\colon X\to X$ and let $\alpha\colon H\to G$ be a group
homomorphism. Then $\ind{\alpha} X:= G\times_H X $ is a
$G$-CW-complex. It is proper if $X$ is proper. The map
$\ind{\alpha}\colon X = H\times_H X \xrightarrow{\alpha \times \id}
G\times_H X = \ind{\alpha} X $ induces a map $\Pi(\ind{\alpha})$ of
the fundamental categories. We obtain a group homomorphism
$\alpha_*:=K_0(\Pi(\ind{\alpha})_*)\colon U^R_H(X,f)\to
U^R_G(\ind{\alpha}X,\ind{\alpha}f)$.

We defined the invariant $(U^R_G(X,f), u^R_G(X,f))$ because it has
certain good properties. It is a functorial equivariant
Lefschetz invariant, the equivariant generalization of a
functorial Lefschetz invariant~\cite[Definition~2.3]{lueck99}.

\begin{Def}\label{def2.4.1}
A \bet{functorial equivariant Lefschetz invariant}\index{functorial
equivariant Lefschetz invariant} on the family of categories $\GCWcat$
of finite proper $G$-CW-complexes for discrete groups $G$ is a pair
$(\Theta,\theta)$ that consists of
\begin{itemize}
\item A family $\Theta$ of functors
\[
 \Theta_G\colon \EndoGCWcat \to \Ab
\]
which is compatible with the induction structure, i.e., for an
inclusion $\alpha\colon G \to K$ there is a group homomorphism
$\Theta(\alpha)\colon \Theta_G(X,f) \to
\Theta_K(\ind{\alpha}X,\ind{\alpha}f)$ for every $(X,f)\in
\EndoGCWcat$. We want the equation
$\Theta(\alpha)\Theta_G(h)=\Theta_K(\ind{\alpha}h)\Theta(\alpha)$ to
hold for any morphism $h\colon(X,f)\to (Y,g)$.
\item A family $\theta$ of functions $\theta_G\colon (X,f)\mapsto
\theta_G(X,f)\in \Theta_G(X,f)$.
\end{itemize}
such that the following holds:
\begin{enumerate}
\item Additivity\\ For a $G$-pushout with $i_2$ a $G$-cofibration
\[
\xymatrix{ (X_0,f_0) \ar[r]^-{i_1} \ar[d]_{i_2} \ar[dr]^{j_0} &
(X_1,f_1) \ar[d]^{j_1} \\ (X_2,f_2) \ar[r]^-{j_2} & (X,f) }
\]
we obtain in $\Theta_G(X,f)$ that
\[
\theta_G(X,f)=\Theta_G(j_1)\theta_G(X_1,f_1) +
\Theta_G(j_2)\theta_G(X_2,f_2) - \Theta_G(j_0)\theta_G(X_0,f_0).
\]
\item $G$-Homotopy invariance\\ If $h_0,h_1\colon (X,f)\to (Y,g)$ are
two $G$-maps that are $G$-homotopic in $\EndoGCWcat$, then
\[
\Theta_G(h_0)=\Theta_G(h_1)\colon \Theta_G(X,f)\to \Theta_G(Y,g).
\]
\item Invariance under $G$-homotopy equivalence\\ If $h\colon (X,f)\to
(Y,g)$ is a morphism in $\EndoGCWcat$ such that $h\colon X\to Y$ is a
$G$-homotopy equivalence, then
\begin{eqnarray*}
\Theta_G(h)\colon \Theta_G(X,f) & \xrightarrow{\cong} &
\Theta_G(Y,g)\\ \theta_G(X,f) & \mapsto & \theta_G(Y,g).
\end{eqnarray*}
\item Normalization: We have $\theta_G(\emptyset,\id_{\emptyset})=0\in
\Theta_G(\emptyset,\id_{\emptyset}).$
\item If $\alpha\colon G\to K$ is an inclusion of groups, then
\[
\alpha_*\theta_G(X,f)=\theta_K(\ind{\alpha}X,\ind{\alpha}f)\in
\Theta_K(\ind{\alpha}X,\ind{\alpha}f).
\]
\end{enumerate}
A natural transformation $\tau\colon (\Theta,\theta) \to (\Xi,\xi)$ of
functorial equivariant Lefschetz invariants is a family of natural
transformations $\tau_G\colon \Theta_G\to \Xi_G$ of functors from
$\EndoGCWcat$ to $\Ab$ for discrete groups $G$ that preserves all
structure.
\end{Def}

\begin{Prop}\label{propufunctequLI}
The invariant $\bigl(U^R_G(X,f),u^R_G(X,f)\bigr)$ is a functorial
equivariant Lefschetz invariant on the family of categories $\GCWcat$
of finite proper $G$-CW-complexes for discrete groups $G$.
\end{Prop}

\begin{proof}
1. Additivity: From~\cite[Lemma~13.7]{lueck89} one knows that for the
corresponding $G$-pushout of $G$-CW-complexes we obtain a based exact
sequence of $R\Pi(G,X)$-chain complexes
\begin{align*}
0\to {j_{0}}_*C^c(\widetilde{X_0}) &
\xrightarrow{j_{1*}C^c(\widetilde{i_1})\oplus
j_{2*}C^c(\widetilde{i_2})} j_{1*}C^c(\widetilde{X_1})\oplus
j_{2*}C^c(\widetilde{X_2}) \\
&\xrightarrow{C^c(\widetilde{j_1})-C^c(\widetilde{j_2})}C^c(\widetilde{X})\to
0,
\end{align*}
where $j_{0*}$, $j_{1*}$ and $j_{2*}$ denote induction with $j_{0}$,
$j_{1}$ and $j_{2}$ respectively. By Lemma~\ref{cd} we obtain a short
exact sequence $0\to {j_0}_* C^c(\widetilde{f_0}) \to
j_{1*}C^c(\widetilde{f_1})\oplus j_{2*}C^c(\widetilde{f_2}) \to
C^c(\widetilde{f}) \to 0$. We conclude using the additive properties
of the algebraic invariant.

2. $G$-Homotopy invariance: Let $H\colon X\times I\to Y$ be a
$G$-equivariant homotopy between $h_0$ and $h_1$. We want to show that
$U^R_G(h_0)=U^R_G(h_1)$. We know that $\Pi(G,h_0)=\Pi(G,H)\circ
\Pi(G,i_0)$ and that $\Pi(G,h_1)=\Pi(G,H)\circ \Pi(G,i_1)$, so it
suffices to prove that $K_0(\Pi(G,i_0)_*(g))=K_0(\Pi(G,i_1)_*(g))$ for
all $g\in \encat{\phi}{\Pi(G,X)}$. There is a natural equivalence
$t\colon \ind{\Pi(G,i_0)}\xrightarrow{\sim}\ind{\Pi(G,i_1)}$ induced
by the paths $v_x$ from $i_1\circ x$ to $i_0\circ x$ given by the map
$x\times (1-\id_I)\colon G/H\times I \to X\times I$ at every object
$x$ of $\Pi(G,X)$. It induces an isomorphism
$\Pi(G,i_0)_*(g)\xrightarrow{\sim}\Pi(G,i_1)_*(g)$ which implies the
desired equality.

3. Invariance under $G$-homotopy equivalence: If $h\colon (X,f)\to
(Y,g)\in \EndoGCWcat$ such that $h\colon X\to Y$ is a $G$-homotopy
equivalence, then $\Pi(G,h)$ is an equivalence of categories. By Proposition~\ref{equivofcat}
 the induced functor $\Pi(G,h)_*\colon
\encat{\phi}{\Pi(G,X)}\to \encat{\psi}{\Pi(G,Y)}$ is an equivalence of
categories. Thus we know that
$U_G(h):=K_0(\Pi(G,h)_*)$ is a bijection. We need to show that
$u_G(X,f)$ maps to $u_G(Y,g)$ under this map.

By Lemma~\ref{cd}, $h$ induces a map $C^c(\widetilde{h})$ from
$\Pi(G,h)_*C^c(\widetilde{f})$ to $C^c(\widetilde{g})$. This map is an
$R\Pi(G,X)$-chain homotopy equivalence. By Lemma~\ref{homequnull}
we conclude $u(\Pi(G,h)_* C^c(\widetilde{f}))=u(
C^c(\widetilde{g}))$, and the claim follows.

4. Normalization: We have $\encat{\id}{\emptyset}=\{\pt\}$ and
$K_0(\{pt\})=\{0\}$.

5. Let $\alpha\colon G\to K$ be an inclusion of groups. We want to
show that $\alpha_* u_G(X,f) = u_K(\ind{\alpha}X,\ind{\alpha}f)$. By
definition we have $\alpha_* u_G(X,f) =
u(\Pi(\ind{\alpha})_*C^c(\widetilde{f}))$ and
$u_K(\ind{\alpha}X,\ind{\alpha}f)=u(C^c(\widetilde{\ind{\alpha}f}))$. There is a natural equivalence $T\colon \ind{\Pi(\ind{\alpha})_*}
\widetilde{X} \to \widetilde{\ind{\alpha}X}$ of
$\Pi(K,\ind{\alpha}X)$-spaces. We check that
$\widetilde{\ind{\alpha}f}T=({\psi}^* T) \Pi(\ind{\alpha})_*
\widetilde{f}$. We can now apply the cellular chain complex functor
$C^c$ to obtain a natural equivalence which maps
$\Pi(\ind{\alpha})_*C^c(\widetilde{f})$ to
$C^c(\widetilde{\ind{\alpha}f})$ and induces the desired equality.
\end{proof}

\section{Universality} \label{sec3}

Now we show that the invariant $(U^\Z_G, u^\Z_G)$ has a universal
initial property among all functorial equivariant Lefschetz
invariants.

\begin{Def}
A functorial equivariant Lefschetz invariant $(U_G,u_G)$ is called
\bet{universal}\index{universal functorial equivariant Lefschetz
invariant} if for any functorial equivariant Lefschetz invariant
$(\Theta_G,\theta_G)$ there is precisely one family of natural
transformations $\tau_{G}\colon U_{G}\to \Theta_{G}$ such that
$\tau_{G (X,f)}\colon U_G(X,f)\to \Theta_G(X,f)$ sends $u_G(X,f)$ to
$\theta_G(X,f)$ for any object $(X,f)$ in $\EndoGCWcat$, for any
discrete group $G$, and such that the equality $\tau_K \circ U(\alpha)
= \Theta(\alpha) \circ \tau_G$ holds for inclusions $\alpha\colon G\to
K$.
\end{Def}

The goal of this section is the proof of Theorem~\ref{universalitytheorem}.

\begin{Thma}
The pair $(U^\Z_G,u^\Z_G)$ is the universal functorial equivariant
Lefschetz invariant on the family of categories $\EndoGCWcat$ for
discrete groups $G$.
\end{Thma}

The proof is in analogy to L\"uck~\cite[Section~4]{lueck99}, of
which it is the equivariant generalization. Before starting, we
introduce notation.

Let $X$ be a $G$-space. A \emph{retractive $G$-space $Y$ over $X$} is
a triple $Y=(Y,i,r)$ which consists of a $G$-space $Y$, a
$G$-cofibration $i\colon X\to Y$ and a $G$-map $r\colon Y\to X$
satisfying $r\circ i=\id_X$. Given a retractive $G$-space $Y$ over
$X$, we define retractive $G$-spaces $Y\times_X [0,1]$ and $C_X Y$ to
be the pushouts
\[
\xymatrix{ X\times [0,1] \ar[r]^-{\pr_X} \ar[d]_{i\times \id} & X
\ar[d] & Y\times \{ 1\} \ar[r]^-{r} \ar[d]_{\incl} & X \ar[d] \\
Y\times [0,1] \ar[r] & Y\times_X [0,1] & Y\times_X [0,1] \ar[r] & C_X
Y, }
\]
where $[0,1]$ is endowed with the trivial $G$-action. The inclusions
of $X$ into $Y\times_X [0,1]$ and into $C_X Y$ are the right vertical
maps, the retractions $Y\times_X [0,1]\to X$ and $\widehat{r}\colon C_X
Y\to X$ are induced by the retraction $r\colon Y\to X$ by the
pushout property.

Define the retractive $G$-space $\Sigma_X Y$ to be the pushout of two
copies of the inclusion $\widehat{i}\colon Y \to C_X Y$ induced by
$Y\times \{0\} \to Y\times [0,1]$. The retraction is induced by
$\widehat{r}$ by the pushout property. The composition
$\widehat{i}\circ i\colon X\to C_X Y$ is a $G$-homotopy equivalence
relative $X$ with the retraction of $C_X Y$ onto $X$ as $G$-homotopy
inverse relative $X$.

Given two retractive spaces $Y$ and $Z$ over $X$ and a
$G$-endomorphism {$f\colon X\to X$,} define $[(C_XY,Y),(C_X Z,Z)]^G_f$
to be the set of $G$-homotopy classes relative $X$ of maps of pairs
$(\widehat{g},g)\colon (C_X Y, Y)\to (C_X Z,Z)$ that induce the given
endomorphism $f$ on $X$, i.e., such that $g i_Y=i_Z f$.

We define the $\Z\Pi(G,X)$-chain complex
$C^c(\widetilde{Y},\widetilde{X})$ by
\[
C^c(\widetilde{Y},\widetilde{X}):=\coker\bigl(C^c(\widetilde{X})\xrightarrow{C^c(\widetilde{i})}
\rrres{\Pi(G,i_Y)}C^c(\widetilde{Y})\bigr).
\]

We call a retractive $G$-space $Y$ over $X$ a \bet{$d$-extension} if
$Y$ is obtained from~$X$ by attaching finitely many cells in dimension
$d$. If $Y$ is a $d$-extension of~$X$ and $d\geq 2$, then we have the
short exact sequence
\[
0 \to C^c(\widetilde{X})\xrightarrow{C^c(\widetilde{i})}
\rrres{\Pi(G,i_Y)}C^c(\widetilde{Y}) \to
C^c(\widetilde{Y},\widetilde{X})\to 0,
\]
where $C^c(\widetilde{X})$ and $\rrres{\Pi(G,i_Y)}C^c(\widetilde{Y})$
are finite free $\Z\Pi(G,X)$-chain complexes equipped with the
cellular equivalence class of bases~\cite[Definition~13.3 and
Example~9.18]{lueck89}. Since the above sequence is based split exact,
$C^c(\widetilde{Y},\widetilde{X})$ is a chain complex concentrated in
degree $d$, and $C^c_d(\widetilde{Y},\widetilde{X})$ is a finite free
$\Z\Pi(G,X)$-module with an equivalence class of bases given by the
cells of $Y\setminus X$~\cite[Example~9.18]{lueck89}.

If $Y$ is a $d$-extension of $X$ and $d\geq 2$, then
$\Pi(G,i_Y)\colon \Pi(G,X)\to \Pi(G,Y)$ is an equivalence of
categories with inverse $\Pi(G,r_Y)$. This implies that
$\rrres{\Pi(G,i_Y)}C^c(\widetilde{Y})$ is a finite free
$\Z\Pi(G,X)$-module.

For any map $g\colon Y\to Z$ of $d$-extensions of $X$ ($d\geq 2$)
such that $g i_Y = i_Z f$ we obtain a lift $\widetilde{g}\colon
\rrres{\Pi(G,i_Y)}\widetilde{Y} \to \rrres{\Pi(G,i_Z)}\widetilde{Z}\circ
\phi$ uniquely determined by the fact that it induces $\widetilde{f}$
on $\widetilde{X}$. We have $C^c(\widetilde{g})
C^c(\widetilde{i_Y})=(\phi^*C^c(\widetilde{i_Z}))C^c(\widetilde{f})$ and
therefore an induced map $C^c(\widetilde{g},\widetilde{f})\colon
C^c(\widetilde{Y},\widetilde{X})\to
C^c(\widetilde{Z},\widetilde{X})\circ \phi$ of $\Pi(G,X)$-chain
complexes. This leads to a bridge between geometry and algebra, in analogy to
L\"uck~\cite[Lemma~4.2]{lueck99}.

\begin{Lem}\label{l2}
Let $Y$ and $Z$ be $d$-extensions of the $G$-space $X$ with $d\geq
2$. Then there is a bijective map
\begin{eqnarray*}
\eta: \bigl[(C_XY,Y),(C_XZ,Z)\bigr]_f^G & \to &
\Mor_{\Z\Pi(G,X)}\bigl(C_d^c(\widetilde{Y},\widetilde{X}),
C^c_d(\widetilde{Z},\widetilde{X})\circ \phi\bigr)\\
{}[(\widehat{g},g)] & \mapsto & C^c_d(\widetilde{g},\widetilde{f}).
\end{eqnarray*}
\end{Lem}

\begin{proof}
Choose a $G$-pushout
\begin{eqnarray} \label{diag1}
\xymatrix{ {}\coprod_{i\in I} G/H_i \times S^{d-1}
\ar[rr]^-{\coprod_{i\in I}q_i} \ar[d] & & X \ar[d]^{i_Y} \\
{}\coprod_{i\in I} G/H_i \times D^{d} \ar[rr]^-{\coprod_{i\in I}Q_i} &
& Y. \\ }
\end{eqnarray}
Define a $G$-map $p_i\colon G/H_i \times S^d \to Y$ by setting
$p_i|_{G/H_i\times S^d_+}:=Q_i$ and $p_i|_{G/H_i\times S^d_-}:=r\circ
Q_i$. Then we can see $C_X Y$ as the pushout
\[
\xymatrix{ {}\coprod_{i\in I} G/H_i \times S^{d}
\ar[rr]^-{\coprod_{i\in I}p_i} \ar[d] & & Y \ar[d]^{\widehat{i_Y}} \\
{}\coprod_{i\in I} G/H_i \times D^{d+1} \ar[rr]^-{\coprod_{i\in I}P_i}
& & C_X Y. \\ }
\]
Analogously to L\"uck~\cite[Lemma~4.2]{lueck99}, we obtain an
isomorphism
\begin{eqnarray*}
\mu\colon \bigl[(C_XY,Y),(C_XZ,Z)\bigr]^G_f & \xrightarrow{\sim} &
\prod_{i\in I}
\pi_d\bigl(Z^{H_i}(f(x_i)),X^{H_i}(f(x_i)),f(x_i)\bigr)\\ &
\xrightarrow{\sim} & \prod_{i\in I}\bigl(
C^c_d(\widetilde{Z},\widetilde{X})\circ \phi\bigr)(x_i)\\ &
\xrightarrow{\sim} &
\Mor_{\Z\Pi(G,X)}\bigl(C^c_d(\widetilde{Y},\widetilde{X}),C^c_d(\widetilde{Z},\widetilde{X})\circ
\phi\bigr).
\end{eqnarray*}
Under this isomorphism,
\begin{eqnarray*}
{}[(\widehat{g},g)] & \mapsto & \bigl([g\circ Q_i(1H_i,-),f\circ
q_i(1H_i,-)]\bigr)_{i\in I} \\ & \mapsto &
\Bigl((\widetilde{g},\widetilde{f})\bigl(\widetilde{Q_i}(1H_i,[D^d]),\widetilde{q_i}(1H_i,[S^d])\bigr)\Bigr)_{i\in
I}=\bigl(C^c_d(\widetilde{g},\widetilde{f})(x_i)\bigr)_{i\in I}\\ &
\mapsto & C^c_d(\widetilde{g},\widetilde{f}).
\end{eqnarray*}
\vskip-4.5ex
\end{proof}

\begin{proof}[Proof of Theorem~\ref{universalitytheorem}]
We have already shown that the pair $(U^\Z_G,u^\Z_G)$ is a functorial
equivariant Lefschetz invariant. It remains to prove universality. For
every functorial equivariant Lefschetz invariant $(\Theta_G,\theta_G)$
we need to find a natural transformation $\tau_G\colon U^\Z_G\to
\Theta_G$ such that $u^R_G(X,f)$ maps to $\theta_G(X,f)$ for all
$(X,f)$ in $\EndoGCWcat$, for discrete groups $G$, and such that
$\tau_K \circ U(\alpha) = \Theta(\alpha) \circ \tau_G$ for inclusions
$\alpha\colon G\to K$.

We define $\tau$ by translating the algebraic data encoded by $U^\Z_G$
back into geometric information using Lemma~\ref{l2}. 
We proceed in eight steps. We omit details which are completely
analogous to~\cite[Section~4]{lueck99}.

{\bf Step 1:} For any $d$-extension $Y$ of $X$, with $d\geq 1$, we
define
\begin{eqnarray*}
\tau_Y \colon \bigl[(C_X Y,Y),(C_X Y,Y)\bigr]^G_f & \to &
\Theta_G(X,f) \\ {}[(\widehat{g},g)] & \mapsto &
\Theta_G(\widehat{i}\circ i)^{-1} \Theta_G(\widehat{i}) \bigl(\theta_G
(g)\bigr).
\end{eqnarray*}
The maps $g$ and $\widehat{g}$ are only defined up to $G$-homotopy
(relative $X$), but because of the $G$-homotopy invariance of
$(\Theta_G,\theta_G)$ this does not play a role.

{\bf Step 2:} We define a map $ \tau_d\colon U^\Z_G(X,f)\to
\Theta_G(X,f) $ for $d\geq 2$. Let $[a\colon M\to M\circ \phi]\in
U^\Z_G(X,f)$. We have $M\cong \bigoplus_{i\in I} \Z\Pi(G,X)(?,x_i)$
with $(x_i\colon G/H_i\to X)\in \Ob\Pi(G,X)$~\cite[p.~167]{lueck89}
since $M$ is a finite free $\Z\Pi(G,X)$-module. Set $q_i:=x_i\circ
\pr_{G/H_i}\colon G/H_i \times S^{d-1} \to X$ and define $Y$ to be the
pushout defined by the maps $q_i$ as in diagram~\ref{diag1}.

Take an isomorphism $c\colon M\xrightarrow{\sim} \bigoplus_{i\in
I} \Z\Pi(G,X)(?,x_i)\xrightarrow{\sim}
C_d^c(\widetilde{Y},\widetilde{X})$. Setting
$[(\widehat{g},g)]:=\eta^{-1}(\phi^* cac^{-1})$ we have $[a]=[\phi^* c
a
c^{-1}]=[\eta([\widehat{g},g])]=[C^c_d(\widetilde{g},\widetilde{f})]=[(-1)^dC^c(\widetilde{g},\widetilde{f})]$. We
define $\tau_d([a]):=(-1)^{d}
\bigl(\tau_Y([(\widehat{g},g)])-\theta_G(f) \bigr)$ and use
$G$-homotopy invariance to show that this is independent of the
choices made.

{\bf Step 3:} We check that $\tau_d$ is compatible with the additivity
relation using the union of $G$-spaces over $X$: For a split short
exact sequence $0\to a_1 \to a \to a_2 \to 0$ we have
$[a_0]+[a_1]=[a]$, and putting $Y:=Y_0\cup_X Y_1$ we can calculate
that $\tau_d([a_0])+\tau_d([a_1])=\tau_d([a])$.

{\bf Step 4:} The map $\tau_d$ is independent of the integer $d\geq
2$. Namely, let $[a]\in U^\Z_G(X,f)$ and let $(Y,g)$ be a
$d$-extension of $X$ such that
$[a]=[C^c_d(\widetilde{g},\widetilde{f})]$. In order to show that
$\tau_d([a])=\tau_{d+1}([a])$, we describe a suspension map
\begin{eqnarray*}
\Sigma_X: \bigl[(C_X Y, Y),(C_X Z,Z)\bigr]^G_f & \to & \bigl[(C_X
\Sigma_X Y, \Sigma_X Y), (C_X \Sigma_X Z, \Sigma_X Z)\bigr]^G_f\\
{}[(\widehat{g},g)] & \mapsto & [(\widehat{\Sigma_X g},\Sigma_X g)].
\end{eqnarray*}
in analogy to L\"uck~\cite[Section~4]{lueck99}, with $\Sigma_X g:=
\widehat{g}\cup_g \widehat{g}$. We have
\begin{eqnarray*}
\tau_d([a]) & = &
(-1)^d\bigl(\tau_Y([(\widehat{g},g)])-\theta_G(f)\bigr)\\
\mbox{and}\qquad\tau_{d+1}([a]) & = & (-1)^{d+1}\bigl(\tau_{\Sigma_X
Y}([(\widehat{\Sigma_X g},\Sigma_X g)])-\theta_G(f)\bigr).
\end{eqnarray*}
We obtain the desired equality by calculating that
\[
\tau_{\Sigma_X Y}\bigl([(\widehat{\Sigma_Xg},\Sigma_X
g)]\bigr)-\theta_G(f) =
-\bigl(\tau_Y([(\widehat{g},g)])-\theta_G(f)\bigr).
\]

{\bf Step 5:} We show that $\tau$ is a natural transformation from
$U^\Z_G$ to $\Theta_G$. Let $h\colon X_1 \to X_2$ be a $G$-equivariant
map between spaces with endomorphisms $f_1$ and $f_2$ respectively
such that $f_2h=hf_1$. Let $[a]\in U_G^\Z(X_1,f_1)$, and let $Y_1$ be
a $d$-extension ($d\geq 2$) of $X_1$ with endomorphism $g_1$ extending
$f_1$, i.e., $g_1|_{X_1}=f_1$, such that
$[C^c_d(\widetilde{g_1},\widetilde{f_1})]=[a]$. Define $(Y_2,g_2)$ as
the pushout of $ (Y_1,g_1)\xleftarrow{i_1} (X_1,f_1)
\xrightarrow{h}(X_2,f_2)$. The map
$C^c_d(\widetilde{g_2},\widetilde{f_2})$ can be used as a
representative of $U^{\Z}_G(h)([a])$. This is shown using the fact
that $\Pi(G,i_1)$ and $\Pi(G,i_2)$ are equivalences of categories,
Proposition~\ref{equivofcat} and additivity.

{\bf Step 6:} The natural transformation $\tau\colon U^\Z_G\to
\Theta_G$ maps $u^\Z_G(X,f)\in U^\Z_G(X,f)$ to $\theta_G(X,f)\in
\Theta_G(X,f)$. For a finite $G$-CW-complex $X$ with endomorphism
$f$, let $Y_n$ be the pushout of $X_n\xleftarrow{k_{n-1}}
X_{n-1} \xrightarrow{j_{n-1}} X$, where all arrows are canonical
inclusions. There is a canonical retraction $r_n\colon Y_n\to X$
induced by the inclusions of $X_{n-1}$, $X_n$ and $X$ into
$X$. Defining $f_k\colon X_k\to X_k$ as the restriction of $f$ to
$X_k$, one obtains an endomorphism $g_n\colon Y_n\to Y_n$ as
$g_n=f_n\cup_{f_{n-1}} f$. Additivity and application of
$\Theta_G(r_n)$ yields
\[
\Theta_G(r_n)\theta_G(g_n)=\Theta_G(i_n)\theta_G(f_n) + \theta_G(f) -
\Theta_G(i_{n-1}) \theta_G(f_{n-1}) .
\]
Summing up, we obtain
\[
\sum_{n=0}^{\dim(X)}
\bigl(\Theta_G(r_n)\theta_G(g_n)-\theta_G(f)\bigr)=\theta_G(f) \in
\Theta_G(X,f).
\]
The analogous equation holds for $(U^\Z_G,u^\Z_G)$. If $n\geq 2$, we
have
\[
{}[C^c(\widetilde{g_n},\widetilde{f})] =
[\rrres{\Pi(G,i_n)}C^c(\widetilde{g_n})]-[C^c(\widetilde{f})] =
U^\Z_G(r_n)u^\Z_G(g_n) - u^\Z_G(f)
\]
because $\Pi(G,i_n)$ is an equivalence of categories. We conclude that
\begin{eqnarray*}
\tau_{G (X,f)}\bigl(U^\Z_G(r_n)u^\Z_G(g_n)-u^\Z_G(f)\bigr) & = &
\tau_{G
(X,f)}\bigl([(-1)^dC^c_d(\widetilde{g_n},\widetilde{f})]\bigr)\\ & =
&\bigl(\Theta_G(\widehat{i_n}i_n)^{-1} \Theta_G(\widehat{i_n})
\theta_G(g_n)-\theta_G(f)\bigr)\\ & = &
\Theta_G(r_n)\theta_G(g_n)-\theta_G(f).
\end{eqnarray*}
The last equation follows because we have strict commutativity $r_n
g_n = f r_n$. Summing up, we obtain $\tau_{G
(X,f)}(u^\Z_G(f))=\theta_G(f)$.

If $n\in \{0,1\}$, we use the suspension $\Sigma_X Y_n$ of $Y_n$ to
get into the range~$d\geq 2$.

{\bf Step 7:} We have to show that $\tau_K \circ U(\alpha) =
\Theta(\alpha) \circ \tau_G$ for any inclusion $\alpha\colon G \to
K$. By definition of the natural transformation $\tau$ and of
$U(\alpha)$ and using the notation established above, this is
equivalent to
\begin{align*}
&(-1)^{d} \bigl(\Theta_K(\widehat{\ind{\alpha}i}\circ
\ind{\alpha}i)^{-1} \Theta_K(\widehat{\ind{\alpha}i}) \bigl(\theta_K
(\ind{\alpha}g)\bigr)-\theta_K(\ind{\alpha}f) \bigr)\\ &=
\Theta(\alpha) \bigl( (-1)^{d} \bigl(\Theta_G(\widehat{i}\circ i)^{-1}
\Theta_G(\widehat{i}) \bigl(\theta_G (g)\bigr)-\theta_G(f) \bigr)
\bigr).
\end{align*}
The retractive space constructions commute with induction. Since
$\Theta$ is a family of functors on $\EndoGCWcat$ for discrete groups
$G$ which is compatible with the induction structure, we have
$\Theta(\alpha) \Theta_G(h)=\Theta_K(\ind{\alpha}h)
\Theta(\alpha)$. So
\begin{align*}
& \Theta(\alpha) \bigl( (-1)^{d} \bigl(\Theta_G(\widehat{i}\circ
i)^{-1} \Theta_G(\widehat{i}) \bigl(\theta_G (g)\bigr)-\theta_G(f)
\bigr) \bigr) \\ & = (-1)^d \bigl(
\Theta_K(\widehat{\ind{\alpha}i}\circ \ind{\alpha}i)^{-1}
\Theta_K(\widehat{\ind{\alpha}i})
\Theta(\alpha)\bigl(\theta_G(g)\bigr) - \Theta(\alpha)\theta_G(f)
\bigr).
\end{align*}
Since $(\Theta,\theta)$ is a functorial equivariant Lefschetz
invariant, the equation
$\Theta(\alpha)(\theta_G(g))=\theta_K(\ind{\alpha}g)$ and the
analogous equation for $f$ hold, so we are finished.

{\bf Step 8:} The natural transformation $\tau$ is uniquely
determined. Any element $[a]$ in $U^\Z_G(X,f)$ can be realized by a
$d$-extension $Y$ for $d\geq 2$ as
\begin{eqnarray*}
[a] & = & [C^c_d(\widetilde{g},\widetilde{f})]\\ & = & (-1)^d
\bigl(U_G^\Z(\widehat{i}\circ i )^{-1} \circ U_G^\Z(\widehat{i})
(u_G^\Z(g)) - u_G^\Z(f)\bigr),
\end{eqnarray*}
and by $\tau_{G (X,f)}(u^\Z_G(f))=\theta_G(f)$ and functoriality it
follows that
\[
\tau_{G (X,f)}\bigl([a]\bigr)=(-1)^d \bigl(\Theta_G(\widehat{i}\circ
i)^{-1}\circ \Theta_G(\widehat{i})(\theta_G(g))-\theta_G(f)\bigr).
\]
\vskip-4.5ex
\end{proof}

\section{Splitting Results} \label{sec4}

In this section, we derive a direct sum decomposition of the abelian
group $U^R_G(X,f)$, making it more accessible to computations.

For
subgroups $H\leq G$, the fixed point sets $X^H:=\{x\in X \; | \; hx=x
\text{ for all }h\in H\}$ and the restrictions $f^H:=f|_{X^H}\colon
X^H\to X^H$ come into play. The Weyl group $WH:= N_GH/H$ acts on
$X^H$. We show in Theorem~\ref{theorem431} that the group
$
U^R_G(X,f)
$
splits up into summands corresponding to the fixed point sets $X^H$
for $(H)\in \consub(G)$, the set of conjugacy classes of subgroups of $G$.

Let $X^{>H}:=\{ x\in X^H\;|\; G_x\neq H\}$, where $G_x$ denotes the
isotropy group of $x$, and $f^{>H}:=
f|_{X^{>H}}$. We also show that the element $u^R_G(X,f)$ maps to the elements given by
the relative maps $(f^H, f^{>H})$. We even have a finer decomposition
corresponding to the orbits of connected components of $X^H\setminus
X^{>H}$ under $f^H$.

We now start analyzing $U^G_R(X,f):=K_0(\encat{\phi}{\Gamma})$ for
$\Gamma=\Pi(G,X)$ and $\phi=\Pi(G,f)$. We restrict ourselves to this
geometric case for simplicity. The results hold for more general
EI-categories $\Gamma$ with endofunctors~$\phi$~\cite{weberdrarb}.

A partial ordering on the objects of any EI-category $\Gamma$
is given by $x\leq y
\Leftrightarrow \Mor_{\Gamma}(x,y)\neq \emptyset$. For an object
$x\colon G/H \to X$ of $\Pi(G,X)$, we define the type $\text{tp}(x)$ to be
$(H)\in\consub(G)$. The partial
ordering on $\Pi(G,X)$ becomes a partial ordering
according to the type of $x$ and the connected component of
$WH\setminus X^H$ which $x$ maps into. On $\consub(G)$, we
obtain the partial ordering given by $(H)\leq (K)\Leftrightarrow
\Mor(G/H,G/K)\neq \emptyset$. The endofunctor $\phi$ respects the partial ordering on $\Pi(G,X)$
since it is given by composition with $f$.

We work with $\phi$-endomorphisms $g\colon M\to M\circ \phi$ of finite
free $R\Gamma$-modules $M$. Every finite free
$R\Gamma$-module is isomorphic to one which is of the form $\bigoplus_{b\in
  B}R\Gamma(?,\beta(b))$ with $B$ a finite set and $\beta\colon B\to \Ob
\Gamma$ a map.

\begin{Lem}\label{splittinglemma1}
Let $g\colon \bigoplus_{b\in
  B}R\Gamma(?,\beta(b)) \to \bigoplus_{b\in
  B}R\Gamma(?,\beta(b))\circ \phi$. Then for $b_0\in B$ we have
\[
g(R\Gamma(?,\beta(b_0)))\subseteq \bigoplus_{b; \beta(b)\geq
  \beta(b_0)} R\Gamma(?,\beta(b))\circ \phi.
\]
\end{Lem}

\begin{proof}
Let $b_0\in B_i$, and let $v\in R\Gamma\bigl(?,\beta(b_0)\bigr)$. The following diagram commutes (because $g$ is a natural transformation):
\[
\xymatrix@C=3pc{
R\Gamma\bigl(?,\beta(b_0)\bigr) \ar[rrr]^-{g_?|_{R\Gamma(?,\beta(b_0))}} &&& \oplus_{{b}\in B} R\Gamma\bigl(\phi(?),\beta({b})\bigr)\\
R\Gamma\bigl(\beta(b_0),\beta(b_0)\bigr) \ar[rrr]^-{g_{\beta(b_0)}|_{R\Gamma(\beta(b_0),\beta(b_0))}} \ar[u]^{v^*} &&& \oplus_{{b}\in B} R\Gamma\bigl(\phi(\beta(b_0)),\beta({b})\bigr)\ar[u]_{\phi^*(v^*)=(\phi(v))^*}.
}
\]
We know that $R\Gamma\bigl(\phi(\beta(b_0)),\beta({b})\bigr)=0$ unless $\phi\bigl(\beta(b_0)\bigr)\leq \beta({b})$. But $b_0\in B_i$ implies $\beta(b_0)\in W_i$ and thus also $\phi\bigl(\beta(b_0)\bigr)\in W_i$, so $\phi\bigl(\beta\bigl(b_0\bigr)\bigr)\leq \beta({b})$ is only possible if $\beta({b})\in W_j$ for $j\geq i$. So $R\Gamma\bigl(\phi\bigl(\beta(b_0)\bigr),\beta\bigl({b}\bigr)\bigr)=0$ if ${b}\not\in B_j$ for all $j\geq i$. The map $\bigl(\phi(v)\bigr)^*$ respects the direct sum decomposition because it is just precomposition with $\phi(v)$, and thus $g_?(v)=\bigl(\phi(v)\bigr)^*\bigl(g_{\beta(b_0)}(\id_{\beta(b_0)})\bigr)$ is nonzero only in components with ${b}\in B_j$ for $j\geq i$. So $g_?\bigl(R\Gamma\bigl(?,\beta(b_0)\bigr)\bigr)\subseteq \bigoplus_{j\geq i}\bigoplus_{{b}\in B_j} R\Gamma\bigl(\phi(?),\beta({b})\bigr)$ for all $b_0\in B_i$, which implies the statement.
\end{proof}

Set $M_H:=\bigoplus_{b \text{ with} \atop \text{tp}(\beta(b))=(H)}
R\Gamma(?,\beta(b))$. Then $M\cong\bigoplus_{(H)\in \consub(G)} M_H$.

\begin{Cor}
$g(M_H)\subseteq \bigoplus_{(K)\geq (H)} M_K\circ \phi$.
\end{Cor}
Analogous results hold for maps $h\colon M\to M'$ of finite free
$R\Gamma$-modules.

We denote by $\Gamma_H$ the full subcategory of $\Gamma$ containing
all objects $x$ of type $(H)$. We set $\phi_H:= \phi|_{\Gamma_H}\colon
\Gamma_H\to \Gamma_H$. The inclusion map $i_H\colon \Gamma_H\to
\Gamma$ is a map of categories with endofunctors, $\phi i_H=i_H
\phi_H$.

We now define functors $S$ and $E$ that induce the desired
splitting.

\begin{Def}
We define the \emph{extension functor} $E_H$ by
\[
E_H:={i_H}_*\colon \encat{\phi_H}{\Gamma_H}\to \encat{\phi}{\Gamma}.
\]
\end{Def}

This generalizes the extension functor
$E_x$~\cite[Definition~9.28]{lueck89}.

We set $S_HM:=\rrres{\Gamma_H}M_H$. A map $g\colon M\to M\circ \phi$
does not change the type, so it induces a map $g_H\colon M_H\to
M_H\circ \phi$. We define $S_H(g):=\rrres{\Gamma_H}g_H\colon S_H M\to
S_H M\circ \phi_H$. For a morphism $h\colon M\to M'$ with $g' h=(\phi^*h) g$ we set
$S_H h:=\rrres{\Gamma_H}h_H\colon S_HM\to S_HM'$. It is easily checked that $S_H (g')
S_H h = (\phi^*_H S_Hh) S_H (g).$

\begin{Def}
We define the \emph{splitting functor} $S_H$ by
\[
S_H\colon \encat{\phi}{\Gamma}\to \encat{\phi_H}{\Gamma_H}, g\mapsto \rrres{\Gamma_H}g_H.
\]
\end{Def}

This is a variation of the splitting functor
$S_x$~\cite[Definition~9.26]{lueck89}, for
objects $x$ of type $(H)$ we have $\rrres{x}S_H=S_x$. The functors $S_H$
and $E_H$ preserve split exact
sequences, so they induce maps on the level of $K$-theory.

We are mostly interested in
$K_0$ since that is where our invariant lives. But we can treat all
higher $K$-groups since $\encat{\phi}{\Gamma}$ is an exact
category. So we let $K$ stand for any
$K_n$, $n\in \Z$, and define the splitting functors for all algebraic $K$-groups
simultaneously. We obtain
\begin{eqnarray*}
K(S_{H})\colon K(\encat{\phi}{\Gamma})  & \to &
K(\encat{\phi_H}{\Gamma_H})\\
K(E_{H})\colon K(\encat{\phi_H}{\Gamma_H})  & \to &
K(\encat{\phi}{\Gamma}).
\end{eqnarray*}

Let another $G$-CW-complex $X'$ with equivariant endomorphism $f'$ be
given and set $\Gamma':= \Pi(G,X')$ and $\psi:=\Pi(G,f')\colon
\Gamma'\to \Gamma'$. A $G$-equivariant map $l\colon X\to X'$
satisfying $f' l = l f$ induces a functor
$L:=\Pi(G,l)\colon \Pi(G,X)\to \Pi(G,X')$ such that $\phi L=L\phi'$ and
which preserves the type. For every $H\leq G$, the functor $L$ induces
a functor ${L_H}_*\colon \encat{\phi_H}{\Gamma_H}\to
\encat{\phi_H'}{\Gamma_H'}$.

\begin{Def}
We define a functor
\begin{eqnarray*}
\Split K\colon \qquad \qquad \EndoGCWcat & \to & \Ab \\
(X,f) & \mapsto & \bigoplus_{(H)\in \consub(G)}
K(\encat{\phi_H}{\Gamma_H})\\
\bigl(l\colon (X,f)\to (X',f')\bigr) & \mapsto & \bigoplus_{(H)\in \consub(G)} K({L_H}_*).
\end{eqnarray*}
\end{Def}
To have compatible notation, we define
$
K\colon \EndoGCWcat  \to  \Ab $ by $
(X,f)  \mapsto K(\encat{\phi}{\Gamma}),
\bigl(l\colon (X,f)\to (X',f')\bigr)  \mapsto   K({L}_*).
$

For finite subsets $I\subseteq \consub(G)$, set $S_I:=\prod_{(H)\in I} K(S_H)$. The map $\colim_{I\subseteq
  \consub(G)}S_I\colon K(\encat{\phi}{\Gamma})\to
\bigoplus_{(H)\in\consub(G)}K(\encat{\phi_H}{\Gamma_H})$ induces a
natural transformation $S\colon K\to \Split K$. (Naturality can be
checked directly, it also follows
from the proof of Theorem~\ref{theorem431}.) The functors $E_H$
combine to form a natural transformation $E\colon \Split K\to
K$. Similar to~\cite[Theorem~9.34]{lueck89}, we prove the
following theorem.

\begin{Thm} \label{theorem431}
We have inverse pairs of natural equivalences $E$ and $S$ between the
functors
\[
K \text{ and } \Split K\colon  \EndoGCWcat \to \Ab.
\]
I.e., if $X$ is a finite proper $G$-CW-complex with equivariant
endomorphism $f\colon X\to X$, and if $\Gamma=\Pi(G,X)$ and
$\phi=\Pi(G,f)$, then
\[
K(\encat{\phi}{\Gamma})\cong \bigoplus_{(H)\in \consub(G)}
K(\encat{\phi_H}{\Gamma_H}),
\]
where the isomorphism is given by $S$ with inverse
$E$ and is natural in~$(X,f)$.
\end{Thm}

\begin{proof}
We have $S\circ E=\colim_{I} S_I \circ \colim_I E_I =
\colim_I(S_I\circ E_I)$ and $E\circ S=\colim_I(E_I\circ S_I)$. Hence
it suffices to show for any finite $I\subseteq \consub(G)$ that
$S_I\circ E_I=\Id$ and $E_I\circ S_I=\Id$ hold.

{\bf 1)} We show that $S_I\circ E_I=\Id$.

We have $S_I\circ E_I =
\bigoplus_{(H)\in I} K(S_H E_H)$. Thus we need to show that $K(S_H
E_H)=\Id$ for all $(H)\in I$. We know that $S_H E_H \colon
\ff R\Gamma_H \to \ff R\Gamma_H$ is naturally equivalent to the
identity~\cite[Lemma~9.31]{lueck89}. Going through the
definitions, it is easily checked that this
natural equivalence extends to a natural equivalence on
$\encat{\phi_H}{\Gamma_H}$.

{\bf 2)} We show that $E_I\circ
S_I = \Id$.

The proof proceeds inductively over the cardinality of
$I$. The beginning $I=\emptyset$ is trivial. In the
induction step, choose $(H)$ maximal in $I$ and write
$I'=I\setminus\{(H)\}$. Since we are
working with finite free modules, we can restrict ourselves to modules of the
form
$
M=\bigoplus_{b\in B} R\Gamma \bigl(?,\beta(b)\bigr).
$
Setting
$M_{I'}:=\bigoplus_{(K)\in I'} M_K$, we have a split short exact
sequence
\[
0 \to M_H \xrightarrow{\inc_H} M \to M_{I'} \to 0.
\]
Since $(H)$ is maximal, by Lemma~\ref{splittinglemma1}
we know that $g\inc_H=\inc_H g_H$. We obtain an induced map
$g_{I'}:=\pr_{M_{I'}\circ\phi}g|_{M_{I'}}$ on $M_{I'}$. So
\[
0\to g_H \to g \to g_{I'} \to 0
\]
 is a short exact sequence in the
sense of Definition~\ref{defaddinv}.

We call $\Gamma_{I'}\subseteq \Gamma$ the
full subcategory of $\Gamma$ with objects in $\bigcup_{(K) \in
  I'} \Gamma_K$ and set $\phi_{I'}:=\phi|_{\Gamma_{I'}}$. The module
$M_{I'}$ can be restricted to the full
subcategory $\Gamma_{I'}$. Analogously, induction with the inclusion
$\Gamma_{I'}\subseteq \Gamma$ allows us to view an $R\Gamma_{I'}$-module $N$ as an
$R\Gamma$-module. We extend these assignments to functors
\begin{align*}
 G & \colon \encat{\phi}{\Gamma}\to \encat{\phi_{I'}}{\Gamma_{I'}},  g\mapsto \rrres{\Gamma_{I'}}g_{I'}, h \mapsto \rrres{\Gamma_{I'}}(\pr_{M'_{I'}}h|_{M_{I'}})\\
 F & :=(\incl_{\Gamma_{I'}})_* \colon  \encat{\phi_{I'}}{\Gamma_{I'}} \to \encat{\phi}{\Gamma}
\end{align*}
We have $FG(g)=g_{I'}$ and $GF=\Id$. We have
$E_H S_H (g)=g_H$.
The above sequence induces a cofibration
sequence of functors on $\encat{\phi}{\Gamma}$
\[
0\to E_{H} S_{H}  \to \Id \to FG \to 0.
\]
By the additivity theorem~\cite[p.~103-106]{quillen73}, we obtain
\[
K(E_{H} )K(S_{H} ) + K(F) K(G) = \Id.
\]

The following diagram commutes on the outside because of this equation and on the top because of the induction hypothesis $E_{I'}\circ S_{I'} = \Id$.
\begin{tiny}
\[
\xymatrix@C=-4pc{
K(\encat{\phi_{I'}}{\Gamma_{I'}})\oplus K(\encat{\phi_{H}}{\Gamma_H})
\ar[rr]^\Id
\ar[dr]^{S_{I'}\oplus \Id }
&
& K(\encat{\phi_{I'}}{\Gamma_{I'}})\oplus K(\encat{\phi_{H}}{\Gamma_H})
\ar[dd]^{K(F)\oplus K(E_{H} )}\\
& \bigl(\oplus_{(K)\in I'}K(\encat{\phi_K}{\Gamma_K})\bigr)\oplus K(\encat{\phi_H}{\Gamma_H})
\ar[ur]^{E_{I'}\oplus \Id } \ar[dr]^{E} & \\
K(\encat{\phi}{\Gamma}) \ar[rr]^\Id
\ar[uu]^{K(G) \oplus K(S_{H} )} \ar[ur]^{S} & & K(\encat{\phi}{\Gamma}).
}
\]
\end{tiny}
The right triangle commutes if $\bigl(K(F)\oplus K(E_{H})\bigr)\circ (E_{I'}\oplus \Id)=E$, i.e., if
$
K(F)\circ\bigl(\oplus_{(K)\in I'} K((E_{I'})_{K})\bigr) \oplus K(E_{H}) = \oplus_{(K) \in I} K(E_{K}).
$
This is obviously true: The $(E_{I'})_{K}$ on the left land in $\encat{\phi_{I'}}{\Gamma_{I'}}$, and the functor $F$ pushes them forward to $\encat{\phi}{\Gamma}$ where the $E_{K}$ on the right hand side land.

The left triangle commutes if $(K(S_{I'})\oplus \Id)\circ \bigl(K(G)\oplus K(S_{H})\bigr)=S$, i.e., if
$
\prod_{(K)\in I'} K((S_{I'})_{K}\circ G)\oplus K(S_{H})=\prod_{(H) \in I}K(S_{H}).
$
This is easily seen: For $(K)\in I'$,
we have
$S_K\circ G(g)=S_K(\rrres{\Gamma_{I'}}(\pr_{M_{I'}\circ
  \phi}g|_{M_{I'}}))=\rrres{\Gamma_K}((\rrres{\Gamma_{I'}}(\pr_{M_{I'}\circ
  \phi}g|_{M_{I'}}))_K)=\rrres{\Gamma_{K}}g_K =S_K(g)$.
So the left triangle also commutes. We conclude that the bottom
triangle commutes and thus that $E_{I}\circ
S_{I}=\Id$. Therefore $E$ and
$S$ are inverse isomorphisms.

It remains to show naturality. An equivariant map
$l\colon(X_1,f_1)\to (X_2,f_2)$ induces functors $L:=\Pi(G,l)_*\colon \encat{\phi_1}{\Gamma_1}\to
\encat{\phi_2}{\Gamma_2}$ and $L_H:=\Pi(G,l|_{X_1^H})_*\colon
\encat{{\phi_1}_H}{{\Gamma_1}_H}\to \encat{{\phi_2}_H}{{\Gamma_2}_H}$
for all $(H)\in \consub(G)$.
We have $SE\cong \Id$
and $E(\oplus_{(H)}L_H)\cong LE$, so $SLE\cong SE(\oplus_{(H)}L_H) \cong \oplus_{(H)} L_H$ and the following diagram
commutes:
\[
\xymatrix@C=6pc{
K(\encat{\phi_1}{\Gamma_1}) \ar[r]^{K(L_*)}  & K(\encat{\phi_2}{\Gamma_2})
\ar[d]^{S} \\
\protect\bigoplus_{(H)} \encat{{\phi_1}_H}{{\Gamma_1}_H} \ar[u]_{E} \ar[r]^{\bigoplus_{(H)}K(L_H)_*} &
\protect\bigoplus_{(H)} \encat{{\phi_2}_H}{{\Gamma_2}_H}.
}
\]
Since $ES\cong\Id$, this implies $SK(L_*)\cong SK(L_*)ES \cong \bigl(\oplus_{(H)}K(L_H)_*\bigr)S$, and so $S$ is a natural
transformation.
\end{proof}

We can split the groups $K_0(\encat{\phi_H}{\Gamma_H})$ up even
further, according to the $WH$-orbits of connected components
$C\subseteq X^H$ and the action of $f$ on these. We have to distinguish 2 cases:

1) A certain iterate of $f$ sends $C$ into $WH\cdot C$.

2) No iterate of $f$ sends $C$ into $WH\cdot C$.

In case 1, we say that $C$ is \emph{recurring}, and we call
$l(C):=\min\{n\geq 1\;| \; f^n(C)\subseteq WH\cdot C\}$ the \emph{length} of $C$. We
denote by $\Gamma_C$ the full subcategory of $\Gamma$ with objects
isomorphic to $f^{i}(x)\colon G/H\to f^{i}(C)$, for $x\colon G/H\to C$
and $0\leq i\leq l-1$. We
call the set of recurring components~$T$.

In case 2, we say that $C$ is \emph{transient}. We call $ht(C):=\min\{n\in
\N\;|\; f^n(C)\subseteq WH\cdot C'
\text{ with }C'\in T\}$ the \emph{height} of $C$. We denote by $\Gamma_C$ the
full subcategory of objects isomorphic to $x\colon G/H\to C$. (This
corresponds to the orbit $WH\cdot C$.) Every component $C$ has finite height since we are dealing with
finite proper $G$-CW-complexes $X$. Recurring components have height $0$.

We choose a set $\C$ of representatives, i.e., of connected components $C\subseteq X^H$ such that
$\Ob(\Gamma)=\amalg_{C\in \C} \Ob(\Gamma_C)$.

\begin{Thm}
There is an isomorphism of abelian groups
\[
K(\encat{\phi_H}{\Gamma_H})\cong \bigoplus_{C\in\C \text{
    recurrent}}K(\encat{\phi_C}{\Gamma_C}) \oplus \bigoplus_{C\in\C \text{
    transient}}K(\ff R \Gamma_C).
\]
\end{Thm}

\begin{proof}
We define a functor
\[
A\colon \bigoplus_{C\in\C \text{
    recurrent}}\encat{\phi_C}{\Gamma_C} \oplus \bigoplus_{C\in\C \text{
    transient}}K\ff R \Gamma_C\to \encat{\phi_H}{\Gamma_H}
\]
 by
induction with the inclusion of the relevant
subcategories, inserting the $0$-map if there is no endomorphism given. Analogously, we define
\[
B\colon \encat{\phi_H}{\Gamma_H} \to \bigoplus_{C\in\C \text{
    recurrent}}\encat{\phi_C}{\Gamma_C} \oplus \bigoplus_{C\in\C \text{
    transient}}K\ff R \Gamma_C
\]
 by the corresponding restriction. We
proceed as in the proof of Theorem~\ref{theorem431} to show that $A$
and $B$ are
equivalences of categories inverse to each other.

The idea is to split off the transient components, starting from the
top. If the largest appearing height is $k$, we call $\Gamma_{k}$ the full subcategory of $\Gamma_H$ consisting of objects
$x$ lying in connected components of height $k$, and we call
$\Gamma_{<k}$ the full subcategory consisting of those with smaller
height. The inclusion of these subcategories induces induction and
restriction maps $\ind{k}$, $\rrres{k}$, $\ind{<k}$ and $\rrres{<k}$. In the induction step,
the decisive point is that we have a split short exact sequence
\[
\xymatrix{
0 \ar[r] & \ind{{<k}}\rrres{{<k}} M \ar[r]
\ar[d]_{(\incl_{{<k}})_* \rrres{{<k}} (g)} & M \ar[r]
\ar[d]^{g} & \ind{{k}}\rrres{{k}} M \ar[r]
\ar[d]^{0} &
0 \\
0 \ar[r] & (\ind{{<k}}\rrres{\ind{{<k}}} M)\circ \phi \ar@{=}[r] & M\circ \phi
\ar[r] & 0 \ar[r] &
0
}
\]
to which we apply the additivity theorem. We leave the details
to the reader.
\end{proof}

Having established this result, we next replace the groupoids
$\Gamma_C$ by corresponding groups. This is analogous to the
transition from the fundamental groupoid of a topological space to its
fundamental group. Assuming $X$
connected, we choose a basepoint~$x\in X$ and look at the fundamental
group with respect to this basepoint.

In order for the endomorphism
$f\colon X\to X$ to induce an endomorphism $\phi\colon  \pi_1(X,x)\to\pi_1(X,x)$,
we also need to choose a path $v$ from $x$ to $f(x)$. We define
$\phi:=c_v\circ \pi_1(f)$, where $c_v$ is the conjugation map
$c_v\colon \pi_1(X,f(x))\to \pi_1(X,x), \gamma \mapsto v\gamma
v^{-1}$. Here composition is written from left to right, as is usual
for composition of paths. Choosing a path $v$ from $x$ to $f(x)$ corresponds to
choosing a morphism $w=(\sigma,[v])\in\Mor(f(x),x)$.

For transient components $C$, we choose an object $x\colon G/H\to C$
in $\Gamma_C$. This is identified with the point
$x(1H)\in C$. We obtain an equivalence of categories
$\Aut(x)\to \Gamma_C$, where the
group $\Aut(x)$ is viewed as a category. The choice of $x$ does not play a role, we can use induction and
restriction with this equivalence of categories to identify the groups
obtained from different choices.

Now we look at recurrent components $C$. If $l:=l(C)$ is the length of
$C$, choosing a point $x\in C$ gives a sequence $f(x)\in f(C),
\ldots, f^l(x)\in f^l(C)$. There is an element $g\in WH$ such that
$gf^l(C)\subseteq C$. We choose a path $v$ from $x$ to $gf^l(C)$. We
know that an element $g\in WH$ uniquely determines a map
$\sigma_g\colon G/H\to G/H, g'H\mapsto g'gH$~\cite[Proposition~1.14]{tomdieck79}. We set
$w=(\sigma_g,[v])\in \Mor(\phi^l(x),x)$. For $\gamma\in
\Aut(\phi^l(x))$, we set $c_w(\gamma):=w\gamma w^{-1}\in \Aut(x)$. We obtain a
group homomorphism
\[
\phi_{x,w}\colon
\Aut(\phi^{l-1}(x))\xrightarrow{\phi|_{\Aut(\phi^{l-1}(x))}}
\Aut(\phi^l(x))\xrightarrow{c_w} \Aut(x).
\]
The collection
$\Phi:=(\phi|_{\Aut(x)},\ldots,\phi|_{\Aut(\phi^{l-2}(x))},\phi_{x,w})$
is an endomorphism of the disjoint union $\amalg_{i=0}^{l-1} \Aut(\phi^{i}(x))$ which on
every component is a group homomorphism to the next.

We define a category $\encat{\Phi}{\amalg\Aut(\phi^{i}(x))}$. An
object is a pair of sequences $((M_i)_{0\leq i \leq l-1},(g_i)_{0\leq
  i \leq l-1})$. Here $M_i$ is a finite free
$R\Aut(\phi^{i}(x))$-module. For $0 \leq i<l-1$, the map $g_i\colon M_i\to
M_{i+1}\circ \phi|_{\Aut(\phi^{i}(x))}$ is a
$\phi|_{\Aut(\phi^i(x))}$-morphism. For $i=l-1$ the map $g_{l-1}\colon M_{l-1}\to
M_0\circ \phi_{x,w}$
is a $\phi_{x,w}$-morphism. A morphism $h\colon (M_i)_{0\leq i \leq l-1} \to
(M'_i)_{0\leq i \leq l-1}$ between these modules is a sequence $(
h_0,\ldots,h_{l-1})$ of maps $h_i\colon M_i\to M'_i$ such that all
resulting diagrams commute.

\begin{Thm}\label{thfinesplit}
There is an equivalence of categories
\[
\encat{\phi_C}{\Gamma_C}\xrightarrow{\sim}
\encat{\Phi}{\amalg\Aut(\phi^{i}(x))}.
\]
\end{Thm}

The proof is quite technical, but not very insightful. We state the
basic idea and refer the reader to \cite{weberdrarb} for details.

\begin{proof}[Idea of proof]
We use induction and restriction with the inclusion of categories
$ \amalg_{i=0}^{l-1}\Aut(\phi^{i}(x))\to \Gamma_C$ to define
$I_{x,w}\colon \encat{\Phi}{\amalg\Aut(\phi^{i}(x))}\to
\encat{\phi_C}{\Gamma_C}$ and $R_{x,w} \colon
\encat{\phi_C}{\Gamma_C}\to
\encat{\Phi}{\amalg\Aut(\phi^{i}(x))}$. We use the canonical
isomorphism $\rrres{c_w^{-1}}\rrres{\phi^l(x)}M_0\cong
\rrres{x}M_0$ in the definition of $R_{x,w}$ and its inverse in
the definition of $I_{x,w}$, alongside with the functors
$\eta\colon \Id\to \rrres{}\ind{}$ and $\varepsilon\colon
\rrres{}\ind{}\to \Id$ as in Section~\ref{sec1}.

We then show that $I_{x,w}$ and $R_{x,w}$ are equivalences of
categories inverse to each other, by using the triangular identities to
show that the appropriate diagrams commute.
\end{proof}

We combine the results of this section in the following statement.

\begin{Thm}\label{thmgeomsplitting}
Let $G$ be a discrete group, let $X$ be a finite proper $G$-CW-complex
and let $f\colon X \to X$ be a $G$-equivariant endomorphism. Let
$\bigl(U^\Z_G(X,f),u^\Z_G(X,f)\bigr)$ be the universal functorial Lefschetz
invariant of $(X,f)$. For all $H\leq G$ and $C\in\C_H$, the set of
representatives of connected components of $X^H$, we
choose $x_C\colon G/H\to C$. If $C$ is recurrent of length $l$, we also
choose an element $g_C\in WH$ such that $f^l(C)\in g_CC$ and path $v_C$ from $x_C$ to
$g_Cf^l(x_C)$.

Then there is an isomorphism
\begin{eqnarray*}
\zeta\colon\quad U^\Z_G(X,f) & \xrightarrow{\sim} & \bigoplus_{(H) \in \consub(G)} \Bigl(
 \bigoplus_{C\in\C_H \atop \text{ recurrent}}
K_0(\encat{\Phi_C}{\amalg\Aut(\phi^{i}(x_C))}) \\
&& \qquad \qquad\qquad \qquad \oplus \bigoplus_{C\in\C_H \atop \text{ transient}} K_0(\ff \Z \Aut(x_C))\Bigr)\\
u_G^\Z(X,f) & \mapsto & \sum_{(H)\in \consub(G)} \sum_{C\in\C_H}
u_G^\Z(X,f)_{x_C},
\end{eqnarray*}
where for $C$ recurrent (leaving out the modules) we have
\begin{align*}
  u^\Z_G(X,f)_{x_C}  =  \Bigl[ & \bigl(
  C^c(\widetilde{f}|_{\widetilde{X}^H(x_C)},
  \widetilde{f}|_{\widetilde{X}^{>H}(x_C)}),\ldots \\
&  C^c(\widetilde{g_C^{-1}
  f}|_{\widetilde{X}^H(f^{l-1}(x_C))},\widetilde{g_C^{-1}
  f}|_{\widetilde{X}^{>H}(f^{l-1}(x_C))}) \bigr) \Bigr]
\end{align*}
and for $C$ transient we have
\[
u_G^R(X,f)_{x_C} =
\Bigl[C^c\bigl(\widetilde{X}^H(x_C),\widetilde{X}^{>H}(x_C)\bigr)\Bigr].
\]
\end{Thm}

\begin{proof}
The existence of the isomorphism $\zeta$ was established in Theorems~\ref{theorem431} and ~\ref{thfinesplit}. It
remains to identify the image of $u^\Z_G(X,f)$.

On the modules, $\zeta$ is given by $\rrres{x}S_{\Gamma_C}=S_x$. There is a natural isomorphism
$
S_x\bigl(C^c(\widetilde{X})\bigr)\cong C^c\bigl(\widetilde{X}^H(x),\widetilde{X}^{>H}(x)\bigr)
$ of $\Z \Aut(x)$-chain complexes
for $(x\colon G/H\to X)\in
\Ob(\Pi(G,X))$~\cite[Lemma~9.32]{lueck89}. This gives the result for
transient $C$.

For recurrent $C$, the modules are
\[
M_i=\rrres{f^i(x_C)} S_{\Gamma_C}\bigl(C^c(\widetilde{X})\bigr)
\cong C^c\bigl(\widetilde{X}^H(f^i(x_C)),\widetilde{X}^{>H}(f^i(x_C))\bigr).
\]
For $0\leq i \leq l-2$, the morphisms are
\[
g_i=\rrres{f^i(x_C)}
S_{\Gamma_C}\bigl(C^c(\widetilde{f})\bigr)=
C^c(\widetilde{f}|_{\widetilde{X}^H(f^i(x_C))},
\widetilde{f}|_{\widetilde{X}^{>H}(f^i(x_C))}).
\]
The canonical isomorphism used in Theorem~\ref{thfinesplit} is
geometrically the one induced by multiplication with $g_C^{-1}$,
\[
C^c(\widetilde{g_C^{-1}})\colon
C^c\bigl(\widetilde{X}^H(f^l(x_C)),\widetilde{X}^{>H}(f^l(x_C))\bigr) \to
C^c\bigl(\widetilde{X}^H(x_C),\widetilde{X}^{>H}(x_C)\bigr).
\]
So $
g_{l-1}=
C^c\bigl(\widetilde{g_C^{-1}f}|_{\widetilde{X}^H(f^{l-1}(x_C))},
\widetilde{g_C^{-1}f}|_{\widetilde{X}^{>H}(f^{l-1}(x_C))} \bigr).
$
\end{proof}

\section{The Generalized Equivariant Lefschetz Invariant} \label{sec5}

In this section we develop a generalized equivariant Lefschetz
invariant $\bigl(\Lambda_G(X,f),\lambda_G(f)\bigr)$ as the image of the universal functorial equivariant Lefschetz invariant
$\bigl(U^\Z_G(X,f),u^\Z_G(X,f)\bigr)$ under a convenient trace map $\tr_{G (X,f)}$. It is an equivariant
analog of the generalized Lefschetz invariant~\cite{reidemeister, wecken2} and a refinement of the equivariant Lefschetz class
\cite[Definition~3.6]{lueck-rosenberg}.

We start by defining the group $\Lambda_G(X,f)$ that will be the target group of $\tr_{G (X,f)}$.
We are interested in fixed point
information, so only in objects $(x\colon G/H\to X)\in \Pi(G,X)$ with
$X^H(f(x))=X^H(x)$. The splitting obtained in
Theorem~\ref{thmgeomsplitting} can be written as
\[
U^\Z_G(X,f)\cong \bigoplus_{\overline{x}\in \Is \Pi(G,X), \atop
  X^H(f(x))=X^H(x)} K_0(\encat{\phi_x,w}{\Aut(x)}) \oplus
\text{other terms}.
\]
We design $\Lambda_G(X,f)$ in the same way. For objects
$x\in \Ob\Pi(G,X)$ with $X^H(f(x))=X^H(x)$, we always find
morphisms $w=(\id,[v])\in \Mor(f(x),x)$, and we restrict our
attention to morphisms of that form.

\begin{Def}
For $x\in \Ob\Pi(G,X)$ with $X^H(f(x))=X^H(x)$ and a morphism
$w=(\id,[v_x])\in \Mor(f(x),x)$, set
\[
\Z \pi_1\bigl(X^H(x),x\bigr)_{\phi'_{x,w}}:=
\Z\pi_1\bigl(X^H(x),x\bigr)/\phi_{x,w}(\gamma)\alpha \gamma^{-1}\sim \alpha,
\]
where $\alpha\in \pi_1(X^H(x),x)$, $\gamma\in \Aut(x)$ and
$\phi_{x,w}(\gamma)=w\phi(\gamma)w^{-1}\in \Aut(x)$.
\end{Def}

We have $\phi_{x,w}(\gamma)\alpha \gamma^{-1}\in \pi_1(X^H(x),x)$
for all $\gamma\in \Aut(x)$ and $\alpha\in \pi_1(X^H(x),x)$ because
the map $\phi_{x,w}$ does not change the $WH_x$-part of the
morphism $\gamma$ and $\pi_1(X^H(x),x)$ is normal.

We can move from one basepoint to another in the following way: For a morphism $(\sigma,
[t])\in \Mor_{\Pi(G,X)}(x(H),y(K))$, where $X^H(f(x))=X^H(x)$ and
$X^K(f(y))=X^K(y)$, we choose morphisms $w_x=(\id,[v_x])\in
\Mor(f(x),x)$ and $w_y=(\id,[v_y])\in \Mor(f(y),y)$ and set
\begin{eqnarray*}
(\sigma,[t])^*_{w_x,w_y}\colon  \Z \pi_1\bigl(X^K(y),y\bigr)_{\phi_{y,w_y}'} & \to & \Z \pi_1\bigl(X^H(x),x\bigr)_{\phi_{x,w_x}'}\\
\alpha & \mapsto & v_x f(t^{-1}) \sigma^*(v_y^{-1}) \sigma^*(\alpha) t.
\end{eqnarray*}

One easily checks that this map is well-defined. The next lemma shows
that the map $(\sigma,[t])^*_{w_x,w_y}$ does
not depend on the choice of $(\sigma,[t])\in \Mor(x,y)$.

\begin{Lem}\label{lemindepofsigmat}
For $\overline{\alpha}\in \pi_1(X^K(y),y)_{\phi_{y,w_y}'}$ we have
\[
(\sigma,[t])^*_{w_x,w_y} (\overline{\alpha}) = (\tau,
[s])^*_{w_x,w_y}(\overline{\alpha})\in \Z \pi_1\bigl(X^H(x),x\bigr)_{\phi_{x,w_x}'}
\]
for all morphisms $(\sigma,[t]), (\tau,[s])\in \Mor(x,y)$.
\end{Lem}

\begin{proof}
For all $\gamma=(\sigma_\gamma,[v_\gamma])\in \Aut(x)$, we have
\begin{eqnarray*}
(\sigma,[t])^*_{w_y,w_x} (\overline{\alpha}) & = & \overline{v_x f(\gamma) v_x^{-1}
(\sigma,[t])^*_{w_x,w_y} (\alpha) \gamma^{-1}} \\
& = & \overline{v_x {\sigma_{\gamma}^{-1}}^* f(v_\gamma t^{-1})
{\sigma_{\gamma}^{-1}}^* \sigma^* (v_y^{-1}) {\sigma_{\gamma}^{-1}}^*
\sigma^* (\alpha) {\sigma_{\gamma}^{-1}}^* (t v_{\gamma}^{-1})}.
\end{eqnarray*}
Setting $\sigma_{\gamma}=\tau^{-1}\sigma$ and $v_{\gamma}=
(\tau^{-1}\sigma)^*(s^{-1}) t$, we obtain
\[
(\sigma,[t])^*_{w_y,w_x} (\overline{\alpha})
 =  \overline{v_x f(s^{-1}) \tau^* (v_y^{-1}) \tau^*(\alpha) s}\\
 =  (\tau, [s])^*_{w_x,w_y}(\overline{\alpha}).
\]
\vskip-4.5ex
\end{proof}

We can use these maps to change the morphism $w$ and the base point
$x$, so the definition is independent of the choices of $x$ and
$w$. We write $(\sigma,[t])^*$ from now on.

\begin{Def}
\[
\Lambda_G(X,f):=\bigoplus_{\overline{x}\in \Is \Pi(G,X),
  \atop X^H(f(x))=X^H(x)} \Z \pi_1(X^H(x),x)_{\phi_{x,w}'}.
\]
\end{Def}

A trace map
from $K_0\bigl(\encatpz{\phi_{x,w}}{\Aut(x)}\bigr)$ to $\Z
\pi_1(X^H(x),x)_{\phi_{x,w}'}$ will now be defined. (The index fp stands for finitely
generated projective modules.) $\Aut(x)$ is a group extension
$
1\to \pi_1\bigl(X^H(x),x\bigr)\to \Aut(x) \to WH_x \to 1,
$
where $WH_x\leq WH$ is the subgroup fixing $X^H(x)$. There is a commutative diagram
\[
\xymatrix{
1\ar[r] & \pi_1\bigl(X^H(x),x\bigr)\ar[r] \ar[d]_{\phi_{x,w}\big|_{\pi_1(X^H(x),x)}} & \Aut(x) \ar[r] \ar[d]^{\phi_{x,w}} & WH_x \ar[r] \ar[d]^\id &  1\\
1\ar[r] & \pi_1\bigl(X^H(x),x\bigr)\ar[r] & \Aut(x) \ar[r] &  WH_x \ar[r] &  1.
}
\]
We combine the trace maps $\tr_{RG}\colon RG\to R, \sum_{g\in G}r_g
\cdot g \mapsto r_1$~\cite[1.1 and~1.2]{lueck-rosenberg}, applied to the $WH_x$-part, and
$\tr_{(\Z\pi,\phi)}\colon \Z\pi\to \Z\pi_\phi, \sum_{\gamma\in \pi}
n_\gamma \cdot \gamma \mapsto \sum_{\gamma \in \pi} n_\gamma \cdot
\overline{\gamma}$~\cite[3.6]{lueck99}, applied to the
$\pi_1(X^H(x),x)$-part.

We formulate the definition of the trace
map independently of the concrete group $\Aut(x)$ which is given to us geometrically.

\begin{Def}\label{defusefultracemap}
Let $\pi$ and $W$ be discrete groups, and let $G$ be a group
extension $1\to \pi \to G \to W \to 1$.\index{discrete group extension} Let an endomorphism
$\phi\colon G\to G$ be given that restricts to an endomorphism
$\phi_\pi:\pi\to \pi$ and becomes trivial when the normal subgroup
$\pi\leq G$ is divided out. For $R$ be a commutative associative ring with unit, define
$
R \pi_{\phi'}:= R\pi/\sim,
$
where $\phi(\gamma)\alpha \gamma \sim \alpha$ for $\alpha\in \pi$ and
$\gamma\in G$.

Let $\encatp{\phi}{G}$ denote
the category of $\phi$-twisted endomorphisms of finitely generated
projective $RG$-modules.
We define the \bet{trace map}\index{trace map}
\[
\tr_{RG}\colon  \Ob(\encatp{\phi}{G}) \to R \pi_{\phi'}
\]
as follows: On $RG$, we set
$
\tr_{RG}:  R G \to R \pi_{\phi'},
\sum_{g\in G} r_g \cdot g  \mapsto  \sum_{g\in \pi} r_g \cdot \overline{g},
$
where $\overline{\,\cdot\,}\colon \pi\to \pi_{\phi'},g\mapsto \overline{g}$
denotes the projection. Given a $\phi$-twisted endomorphism $u\colon
P\to \rrres{\phi}P$ of a finitely generated projective $RG$-module, we
choose a finitely generated projective $RG$-module $Q$ and an
isomorphism~$v\colon P\oplus Q \xrightarrow{\sim} \bigoplus_{i\in I}
RG$ for a finite indexing set $I$. Then we have a $\phi$-twisted endomorphism
$
\phi^*(v)\circ (u\oplus 0) \circ v^{-1} \colon \bigoplus_{i\in I} RG \to \rrres{\phi}\Bigl(\bigoplus_{i\in I} RG\Bigr),
$
to which a matrix $A=(a_{ij})_{i,j\in I}$ is associated. We define
\[
\tr_{RG}(u):= \sum_{i\in I} \tr_{RG}(a_{ii})\in R\pi_{\phi'}.
\]
\end{Def}

Note that $\alpha\in \pi$
implies that $\phi(\gamma)\alpha \gamma^{-1}\in \pi$ since
$
\pr_W(\phi(\gamma)\alpha\gamma^{-1})=\id_W(\pr_W (\gamma)) \pr_W (\alpha)
\pr_W (\gamma^{-1}) = 1_W$, so $R \pi_{\phi'}$ is well defined. As
usual, the definition of the trace map is independent of the choices involved. This
trace map has properties generalizing~\cite[Lemma~1.3]{lueck-rosenberg}:

\begin{Lem}\label{1.3*}
Let $G$ be a discrete group extension $1\to \pi \to G \to W \to 1$ with
endomorphism $\phi\colon G\to G$ which restricts to $\pi$ and such
that $\phi_W=\id_W$.
\begin{enumerate}
\item \label{1.3*.1}
Let $u\colon P\to Q$ and $v\colon Q\to \rrres{\phi} P$ be $RG$-maps of finitely generated projective $RG$-modules. Then
$
\tr_{RG}(v\circ u)=\tr_{RG}\bigl(\rrres{\phi}(u)\circ v\bigr).
$
\item \label{1.3*.2}
Let $P_1$, $P_2$ be finitely generated projective $RG$-modules. Given
a $\phi$-twisted endomorphism
$
A=\left(  \begin{array}{cc}
u_{11} & u_{12} \\
u_{21} & u_{22} \end{array}  \right)
: P_1 \oplus P_2 \to P_1\oplus P_2
$ we have
$
\tr_{RG}(A)=\tr_{RG}(u_{11})+\tr_{RG}(u_{22}).
$
\item\label{1.3*.3}
Let $u_1, u_2\colon P\to \rrres{\phi}P$ be $\phi$-twisted
endomorphisms of a finitely generated projective $RG$-module $P$
and let $r_1, r_2\in R$. Then
\[
\tr_{RG}(r_1\cdot u_1 + r_2 \cdot u_2)=r_1 \tr_{RG}(u_1) + r_2 \tr_{RG} (u_2).
\]
\item\label{1.3*.4}
Let $\alpha\colon G\to
K$ be a homomorphism of discrete group extensions with endomorphisms as in
Definition~\ref{defusefultracemap} with $\alpha_W$ injective, lying in a commutative diagram
\[
\xymatrix@R=1pc{
1 \ar[r] & \pi \ar[r] \ar[d]^{\alpha_\pi} & G \ar[r] \ar[d]^{\alpha} &
W \ar[r] \ar[d]^{\alpha_W} & 1  \\
1 \ar[r] & K_1 \ar[r] & K \ar[r] & K_2 \ar[r] & 1  .
}
\]
Let
$u\colon P\to \rrres{\phi_G}P$ be a $\phi_G$-twisted endomorphism of a
finitely generated projective $RG$-module $P$. Then induction with
$\alpha$ yields a $\phi_K$-twisted endomorphism $\alpha_* u$ of a
finitely generated projective $RK$-module and
\[
\tr_{RK}(\alpha_* u)=\alpha'_* \tr_{RG}(u),
\]
where $\alpha'_*\colon R\pi_{\phi_G'}\to R(K_1)_{\phi_K'}$ is induced by $\alpha_{\pi}$.
\item\label{1.3*.5}
Let $\alpha\colon H\to G$ be an inclusion of discrete group extensions with endomorphisms as in
Definition~\ref{defusefultracemap} with finite index $[G:H]$, lying in a commutative diagram
\[
\xymatrix@R=1pc{
1 \ar[r] & \pi \ar[r] \ar[d]^{\id} & H \ar[r] \ar[d]^{\alpha} &
H_2 \ar[r] \ar[d]^{\alpha_{H_2}} & 1  \\
1 \ar[r] & \pi \ar[r] & G \ar[r] & W \ar[r] & 1 .
}
\]
Let $u\colon P\to \rrres{\phi_G}P$ be a $\phi_G$-twisted endomorphism of a
finitely generated projective $RG$-module $P$. Then the restriction to
$RH$ with $\alpha$ yields a $\phi_H$-twisted endomorphism $\alpha^* u$
of a finitely generated projective $RH$-module and
\[
\id_* \tr_{RH}(\alpha^* u)= [G:H] \cdot \tr_{RG}(u),
\]
where $\id_*\colon R\pi_{\phi_H'}\to R\pi_{\phi_G'}$ denotes the projection.
\item\label{1.3*.6}
Let the subgroup $H\leq G$ be finite such that $|H|$ is invertible in $R$. Let
$u\colon R[G/H]\to \rrres{\phi} R[G/H]$ be a $\phi$-twisted
endomorphism that sends $1H$ to $\sum_{gH\in G/H} r_{gH}\cdot
gH$. Then $R[G/H]$ is a finitely generated projective $RG$-module and
$
\tr_{RG}(u) = |H|^{-1} \sum_{g\in \pi} r_{gH}\cdot \overline{g} \in R\pi_{\phi'}.
$
In particular,
$
\tr_{RG}\bigl(\id_{R[G/H]} \bigr)  =  |H|^{-1} \cdot \Bigl(\sum_{g\in \pi} 1\cdot \overline{g}\Bigr)\in R\pi_{\phi'}.
$
\end{enumerate}
\end{Lem}

\begin{proof}
Check (1) by calculation, using the fact that for $g,g'\in G$ we have
$
g g'\in U
\Leftrightarrow \phi(g') g\in U.
$
Assertions (2) and (3) are clear by definition, and (4) and (5) are again checked by calculation.
In (6), we are given
$
u\colon  1H \mapsto \sum_{gH\in G/H} r_{gH} \cdot gH = |H|^{-1} \Bigl(\sum_{g\in G} r_{gH} \cdot g\Bigr)H.
$
So
$
\tr_{RG}(u)=|H|^{-1} \tr_{RG}\Bigl(\sum_{g\in G} r_{gH} \cdot g\Bigr)= |H|^{-1}
\sum_{g\in \pi} r_{gH} \cdot \overline{g}. \hfill \qedhere
$
\end{proof}

By assertions~\ref{1.3*.1} and~\ref{1.3*.2} of Lemma~\ref{1.3*}, the
trace map $\tr_{RG}$ is compatible with the relations defining
$K_0(\encatp{\phi}{G})$. So we can use its value on any representative
and define\index{trace map}
\[
\tr_{RG}\colon  K_0(\encatp{\phi}{G})  \to  R \pi_{\phi'},
{}[u] \mapsto  \tr_{RG}(u).
\]

\begin{Rem}\label{remhochschild}
The trace map $\tr_{RG}$\index{trace map} can be seen as a variation of a trace map between $K$-theory and
Hochschild homology with coefficients in the bimodule
$M_\phi=RG(\phi(?),??)$. There is
a trace map $\tr_{(K\to HH)}\colon K(RG; M_\phi)\to HH(RG;M_\phi)$. One
has $HH_0(RG;M_\phi)\cong RG_\phi$, where we define
$RG_\phi:=RG/\phi(\gamma)\beta \gamma^{-1} \sim \beta $ for $\gamma,\beta\in
G$. We have an inclusion $R\pi\to RG$ of group rings. This inclusion is
respected by the $G$-action given by twisted conjugation since
$\beta\in \pi$ implies $\pr_W(\phi(\gamma)\beta
\gamma^{-1})=\pr_W(\gamma)\cdot 1 \cdot \pr_W(\gamma^{-1})=1$, so
$\phi(\gamma)\beta \gamma^{-1}\in \pi$ for all $\gamma\in G$. It
induces an inclusion $R\pi_{\phi'}\to RG_{\phi}$ as a direct
summand. Denoting the restriction to this summand by $r_\pi$, one can
check that $\tr_{RG}=r_\pi \circ {\tr_{(K\to HH)}}_0$.
\end{Rem}

We could now return to our original setting, but we keep the
more general formulation to make a few
observations that will be useful later on.

Let $G$ be a discrete group extension with endomorphism $\phi$ and let $(X,A)$ be a finite proper relative
$G$-CW-complex. Let $R$ be a commutative ring such that the order of
the isotropy group $|G_x|$ is invertible in $R$ for every $x\in
X\setminus A$. Then the cellular $RG$-chain-complex $C^c(X,A)$ is
finite projective. Let $(f,f_0)\colon (X,A)\to (X,A)$ be a $\phi$-twisted
cellular endomorphism. This induces $C^c(f,f_0)\colon C^c(X,A)\to \rrres{\phi} C^c(X,A)$, a $\phi$-twisted
endomorphism of the
cellular chain complex.

\begin{Def}\label{1.4*}
With notation as above, we define the \bet{refined equivariant
  Lefschetz number}\index{refined equivariant Lefschetz number} of $(f,f_0)$ to be
\[
L^{R G}(f,f_0):= \sum_{p\geq 0} (-1)^p \tr_{R G} \bigl(C^c_p(f,f_0)\bigr) \in R \pi_{\phi'}.
\]
\end{Def}
This generalizes the orbifold
Lefschetz number~\cite[Definition~1.4]{lueck-rosenberg}. Writing $\bigl[C^c(f,f_0)\bigr]:=\sum_{p\geq 0} (-1)^p
\bigl[C^c_p(f,f_0)\bigr]\in K_0(\encatp{\phi}{G})$, this definition becomes $
L^{R G}(f,f_0):= \tr_{R G} \bigl([C^c(f,f_0)]\bigr) \in R\pi_{\phi'}.
$

We also have a
refinement of the incidence
number~\cite[1.8]{lueck-rosenberg}.

\begin{Def}
Let $G$ be a discrete group extension with endomorphism $\phi$, let $(X,A)$ be a finite proper
relative $G$-CW-complex, let $e\in I_p(X,A)$ be a $p$-cell and let
$(f,f_0)\colon (X,A)\to (X,A)$ be a $\phi$-twisted cellular
endomorphism. We define the \bet{refined incidence
  number}\index{refined incidence number}
$\inc_\phi(f,e)\in \Z \pi_{\phi'}$ for a $p$-cell $e\in I_p(X,A)$ to
be the ``degree'' of the composition
\begin{align*}
\overline{e}/\partial e & \xrightarrow{i_e} \hspace{-1em}\bigvee_{e'\in I_p(X,A)}
\overline{e'}/\partial e' \xrightarrow{h \sim} X_p/X_{p-1}
\xrightarrow{f} X_p/X_{p-1} \xrightarrow{h^{-1} \sim} \hspace{-1em} \bigvee_{e'\in
  I_p(X,A)} \overline{e'}/\partial e'\\
& \xrightarrow{\pr_{\pi\cdot \overline{e}/\partial e }} \pi \cdot
\overline{e}/\partial e \xrightarrow{\overline{\, \cdot \,}}\pi_{\phi'} \cdot
\overline{e}/\partial e.
\end{align*}
\end{Def}

We have $\inc_{\phi}(f,e)=\inc_{\phi}(f,ge)$ for all $g\in G$, this is ensured by using~$R\pi_{\phi'}$. We have the equation
$
\inc_\phi(f,e)=\sum_{\alpha\in \pi}\inc(\alpha^{-1}f,e)\cdot \overline{\alpha},
$
where the incidence number~\cite[1.8]{lueck-rosenberg} appears on the right hand side.

Let $e\in I_p(X,A)$ be a $p$-cell. Then $C^c(X,A)_e$ is the chain
complex concentrated in degree $p$ with $C^c_p(X,A)_e= R [G/G_e]$. If
\[
C^c_p(f,f_0)|_{e}=\sum_{gG_e\in G/G_e}r_{gG_e}\cdot gG_e\colon R [G/G_e]
\to \rrres{\phi} R [G/G_e] ,
\]
then $\inc_\phi(f,e) = \sum_{g\in \pi} r_{gG_e} \cdot \overline{g}$,
and by assertion \ref{1.3*.5} of Lemma~\ref{1.3*} we have
\[
\tr_{R G}\bigl(C^c_p(f,f_0)|_{e}\bigr)=|G_e|^{-1} \sum_{g\in \pi} r_{gG_e} \cdot \overline{g} = |G_e|^{-1} \inc_\phi(f,e).
\]

This observation helps us prove the following result.

\begin{Lem}\label{1.9*}
Let $G$ be a discrete group extension with endomorphism $\phi$. Let $(X,A)$ be a finite proper
relative $G$-CW-complex. Let $R$ be a commutative ring such that the order of
the isotropy group $|G_x|$ is invertible in $R$ for every $x\in
X\setminus A$. Let $(f,f_0)\colon (X,A)\to (X,A)$ be a
$\phi$-twisted cellular endomorphism. Then
\[
L^{R G}(f,f_0)= \sum_{p\geq 0} (-1)^p \sum_{G\cdot e\in G \setminus I_p(X,A)} |G_e|^{-1} \cdot \inc_{\phi}(f,e) \in R \pi_{\phi'}.
\]
\end{Lem}

\begin{proof}
We calculate
\begin{eqnarray*}
L^{R G}(f,f_0) & = & \sum_{p\geq 0} (-1)^p \sum_{G\cdot e\in G \setminus I_p(X,A)}\tr_{R G}\bigl(C^c_p(f,f_0)|_{e}\bigr)\\
 & = & \sum_{p\geq 0} (-1)^p \sum_{G\cdot e\in G \setminus
   I_p(X,A)}|G_e|^{-1} \cdot \inc_{\phi}(f,e).
\end{eqnarray*}
\vskip-6ex
\end{proof}

We now return to the original setting. For any $x$ in $\Pi(G,X)$ for which 
$X^H(f(x))=X^H(x)$, the group $\Aut(x)$
is a discrete group extension 
$
1\to \pi_1(X^H(x),x) \to \Aut(x) \to WH_x  \to
1
$
 with endomorphism $\phi_{x,w}$ which restricts to $\pi_1(X^H(x),x)$ as
$(\phi_{x,w})_{\pi_1(X^H(x),x)}=c_{w}\circ \pi_1(f^H(x),x)$ and such that
$(\phi_{x,w})_W=\id_{WH_x}\colon WH_x \to WH_x$. We set
\[
{\tr_{G (X,f)}}_x := \tr_{\Z \Aut(x)}\colon
K_0(\encatpz{\phi_{x,w}}{\Aut(x)}) \to \Z
\pi_1\bigl(X^H(x),x\bigr)_{\phi'_{x,w}}.
\]
This definition is independent of the choice of the
element $x\in\overline{x}$ and the morphism $w$. The reason is that
on the domain as well as on the target space we use morphisms to pass
from one choice to another and that these constructions are compatible with
the trace map.

\begin{Def}
Let $G$ be a discrete group, $X$ a finite proper $G$-CW-complex and
$f\colon X \to X$ a $G$-equivariant endomorphism. We
define\index{trace map!generalized}
\[
\tr_{G (X,f)}:=
\bigoplus_{\overline{x}\in \Is\Pi(G,X),\atop X^H(f(x))=X^H(x)} {\tr_{G (X,f)}}_{x} \oplus 0\colon  U^{\Z}_G(X,f) \to \Lambda_G(X,f).
\]
\end{Def}

The groups $\Lambda_G(X,f)$ combine to form a family of functors we want.

\begin{Lem}\label{lemLambdafunctor}
The groups $\Lambda_G(X,f)$ are naturally endowed with the structure of a
family of functors $\Lambda_G$ from $\EndoGCWcat$ to $\Ab$, for discrete groups $G$, which is compatible
with the induction structure.
\end{Lem}

\begin{proof}
Let $l\colon (X,f)\to (Y,g)$ be a $G$-equivariant map between
$G$-CW-complexes with endomorphisms. The map $l$ induces a
functor $\Pi(G,l)$ and a group
homomorphism $\ell_x:=\Pi(G,l)|_{\Aut(x)}\colon \Aut(x)\to
\Aut(l(x))$ for every $x\in \Pi(G,X)$. We have a map
$\ell_x \big|_{\pi_1\bigl(X^H(x),x\bigr)}=\pi_1(l|_{X^H(x)},x)$ which
induces a map $(\ell_x')_*\colon \Z \pi_1(X^H(x),x)_{\phi_x'} \to \Z
\pi_1(Y^H(l(x)),l(x))_{\phi_{l(x)}'}$. We set
\begin{eqnarray*}
\Lambda_G(l)_x := (\ell_x')_* \colon \Z \pi_1\bigl(X^H(x),x\bigr)_{\phi_x'} &
\to & \Z \pi_1\bigl(Y^H(l(x)),l(x)\bigr)_{\phi_{l(x)}'}\\
\sum_{\overline{g}\in \pi_1(X^H(x),x)_{\phi_x'}} r_g \cdot \overline{g} & \mapsto &
\hspace{-2ex}\sum_{\overline{g}\in \pi_1(X^H(x),x)_{\phi_x'}} \hspace{-1em} r_g \cdot \overline{\pi_1(l|_{X^H(x)},x)(g)} .
\end{eqnarray*}
These maps combine to define $\Lambda_G(l)\colon
\Lambda_G(X,f)\to \Lambda_G(Y,g).$

Given an inclusion $\alpha\colon  G\to K$,
the map $\alpha|_{\pi_1(X^H(x),x)_{\phi_G'}}$ induces a map
$\alpha_*'\colon \Z\pi_1(X^H(x),x)_{\phi_G'}\to \Z
\pi_1((\ind{\alpha}X)^H(\ind{\alpha}x),\ind{\alpha}x)_{\phi_K'}$. We
set $\Lambda(\alpha)_x=\alpha'_*$ and obtain $\alpha_*=\Lambda(\alpha)\colon \Lambda_G(X,f) \to
\Lambda_K(\ind{\alpha}X,\ind{\alpha}f)$. So $\Lambda$ is also compatible with the induction
structure.
\end{proof}

The maps $\tr_{G (X,f)}$ respect all structure.

\begin{Prop}\label{proptrnattr}
The collection of maps $\tr_{G (X,f)}$\index{trace map!generalized} is a natural
transformation of families of functors from $\EndoGCWcat$ to $\Ab$, for
discrete groups $G$.
\end{Prop}

\begin{proof}
Let $l\colon (X,f)\to (Y,g)$ be a $G$-equivariant map between
$G$-CW-complexes with endomorphisms. Then $\ell_x:=\Pi(G,l)|_{\Aut(x)}\colon \Aut(x)\to
\Aut(a(x))$ lies in the commutative diagram
\[
\xymatrix{
1 \ar[r] & \pi_1\bigl(X^H(x),x\bigr) \ar[r] \ar[d]^{\pi_1(l|_{X^H(x)},x)} &
\Aut(x) \ar[r] \ar[d]^{\ell_x} & WH_x \ar[r] \ar[d]^{\overline{\ell_x}=\incl} & 1 \\
1 \ar[r] & \pi_1\bigl(Y^H(l(x)),l(x)\bigr) \ar[r] & \Aut\bigl(l(x)\bigr) \ar[r] & WH_{l(x)} \ar[r] & 1 .
}
\]
The map $\overline{\ell_x}\colon WH_x\to WH_{l(x)}$ is an inclusion since the elements in
$WH$ which fix the connected component $X^H(x)$ also fix the connected
component $Y^H(l(x))$, by equivariance and continuity of~$l$. We apply Lemma~\ref{1.3*}, assertion~\ref{1.3*.4}
to obtain for all $u\in\encatpz{\phi_{x,w}}{\Aut(x)}$
\[
\tr_{\Z \Aut(l(x))}(\ell_x)_* (u) =(\ell'_x)_* \tr_{\Z \Aut(x)}(u).
\]

Taking all induction maps
$K_0((\ell_{\overline{x}})_*) \colon  U^\Z_G(X,f)_{\overline{x}}\to
U^\Z_G(Y,g)_{\overline{l(x)}}$ together gives $U^\Z_G(l)=K_0(\Pi(G,l)_*)\colon
U^\Z_G(X,f) \to U^\Z_G(Y,g)$. Combining the above equation for all
$\overline{x}\in\Is\Pi(G,X)$, we arrive at
\[
\tr_{G (Y,g)}\circ U^\Z_G(l) (u)=\Lambda_G(l) \circ \tr_{G (X,f)} (u)
\]
for all $u\in U^\Z_G(X,f)$. So the trace map $\tr_{G (X,f)}$ is a natural transformation of
functors from $\EndoGCWcat$ to $\Ab$.

It remains to show compatibility with the induction structure. Given an inclusion
$\alpha\colon  G\to K$, the functor
$\Pi(\ind{\alpha})\colon \Pi(G,X)\to \Pi(K,\ind{\alpha}X)$ induces a
group
homomorphism $\alpha_*\colon U^\Z_G(X,f)\to
U^\Z_K(\ind{\alpha}X,\ind{\alpha}f)$. Let $x\colon G/H\to X$ be given with $X^H(f(x))=X^H(x)$. We set $W_G H:=
N_G H/H$ and $W_K H=N_K H/H$. The map $\alpha\colon G\to K$ is injective, and so $W_G H \to
W_K H$ is injective. If an element of $W_G H$ fixes the component
$X^H(x)$, then its image fixes $(\ind{\alpha} X)^H
(\ind{\alpha}x)$. This implies that $\alpha|_{N_G H}: N_G H \to N_K H$ is
injective. So $\Pi(\ind{\alpha})|_{W_G H_x}\colon W_G H_x\to W_K
H_{\ind{\alpha}(x)}$, induced by
$\Pi(\ind{\alpha})|_{\Aut(x)}$, is
also injective.

We apply assertion~\ref{1.3*.4} of Lemma~\ref{1.3*}
to obtain $\tr_{\Z K} \alpha_* u = \alpha_*'
\tr_{\Z G} u$ for~$u\in \encatp{\phi_G}{G}$, where the map $\alpha|_{\pi_1(X^H(x),x)_{\phi_G'}}$ induces the homomorphism $\alpha_*'\colon
\Z\pi_1(X^H(x),x)_{\phi_G'}\to \Z
\pi_1((\ind{\alpha}X)^H(\ind{\alpha}x),\ind{\alpha}x)_{\phi_K'}$. Since
$\Lambda(\alpha)_x$ is defined to be $\alpha'_*$, these combine to form the map
$\alpha_*\colon \Lambda_G(X,f) \to
\Lambda_K(\ind{\alpha}X,\ind{\alpha}f)$ such that the desired equation
$\tr_{K (\ind{\alpha}X,\ind{\alpha}f)} \alpha_* = \alpha_*
\tr_{G (X,f)}$ holds on $U^\Z_G(X,f)$.
\end{proof}

Now we define the invariant which contains the fixed point
information we are interested in.
\begin{Def}\label{3.6*}
Let $G$ be a discrete group, let $X$ be a finite proper
$G$-CW-complex, and let $f\colon X\to X$ be a $G$-equivariant
cellular endomorphism. We define the \bet{generalized equivariant
  Lefschetz invariant}\index{generalized equivariant Lefschetz
  invariant} of~$f$ by
\[
\lambda_G(f):=\tr_{G (X,f)}\bigl(u_G^{\Z}(X,f)\bigr)\in \Lambda_G(X,f).
\]
\end{Def}

By Proposition~\ref{proptrnattr}, the collection of the $\tr_{G (X,f)}$ is a natural
transformation from $(U^\Z,u^\Z)$ to $(\Lambda,\lambda)$. The pair
$(\Lambda,\lambda)$ inherits all structure from $(U^\Z,u^\Z)$: It is
also a functorial equivariant Lefschetz invariant on the family of
categories $\GCWcat$ for discrete groups $G$.

\begin{Thm}
The pair $(\Lambda,\lambda)$ is a functorial equivariant Lefschetz
invariant on the family of categories $\GCWcat$ for discrete groups
$G$.
\end{Thm}

\begin{proof}
The natural transformation $\tr_{G (X,f)} \colon U^\Z_G(X,f)\to
\Lambda_G(X,f)$ maps $u^\Z_G(X,f)$ to $\lambda_G(X,f)$. So
$\bigl(\Lambda_G(X,f),\lambda_G(X,f)\bigr)$ is a functorial
equivariant Lefschetz invariant.
\end{proof}

It therefore has all properties stated in
Definition~\ref{def2.4.1}. Since we can define the universal
functorial equivariant Lefschetz invariant for any $G$-equivariant
continuous map $f\colon X\to X$, the same is true for the generalized
equivariant Lefschetz invariant.

We can describe the invariant $\lambda_G(f)$ in a more concrete
way. We see that
$
\tr_{G (X,f)}\bigl(u_G^{\Z}(X,f)\bigr)
= \sum_{\overline{x}\in \Is\Pi(G,X), \atop X^H(f(x))=X^H(x)} L^{\Z \Aut(x)}\bigl(\widetilde{f^H(x)},
                               \widetilde{f^{>H}(x)}\bigr).
$
We can use $\Z$ because $\Aut(x)$ operates freely on
$\widetilde{X^H(x)}\setminus \widetilde{X^{>H}(x)}$.

\begin{Rem}
One can obtain any functorial equivariant Lefschetz invariant by applying a
suitable natural transformation to $(U^\Z_G,u^\Z_G)$. For example, the
equivariant analog of the Lefschetz number is the equivariant
Lefschetz class~\cite[Definition~3.6]{lueck-rosenberg}. The natural
transformation mapping the universal functorial equivariant Lefschetz
invariant to the equivariant Lefschetz class is given by the trace map
$\tr_{G(X,f)}$ followed by an augmentation map $s_{G (X,f)}$ induced by the projection
$\pi_1(X^H(x),x)_{\phi_x'}\to \{1\}$, just like in the
non-equivariant case.
\end{Rem}

\section{The Refined Equivariant Lefschetz Fixed Point Theorem}
\label{sec6} 

In this section, we introduce the generalized local equivariant
Lefschetz class $\lambda_G^{loc}(f)$ in
terms of fixed point data. It is a refinement of the local
equivariant Lefschetz class~\cite[Definition~4.6]{lueck-rosenberg}. We prove
Theorem~\ref{0.2*} which states that $\lambda_G(f)=\lambda^{loc}_G(f)$
under quite general conditions. So $\lambda^{loc}_G(f)$ gives a
concrete geometric description of the fixed point information contained in~$\lambda_G(f)$.

We
briefly assemble the necessary
notation. For a space $X$ with action of a
discrete group $G$, we set
$U^G(X):=\bigoplus_{\overline{x}\in\Is\Pi(G,X)} \Z$. Let $K$ be a finite group. The \bet{Burnside ring} $A(K)$ of $K$ is defined
to be the
Grothendieck ring of finite $K$-sets $S$ with the additive structure
induced by disjoint union and the multiplicative structure induced by the
Cartesian product. Additively, $A(K)=U^K(\pt):=\bigoplus_{(H)\in
  \consub(K)} \Z$.

Let $Z$ be a finite $K$-CW-complex and let
$\psi\colon Z \to Z$ be a $K$-equivariant endomorphism. Then the
\bet{equivariant Lefschetz class with values in the Burnside ring} of $\psi$ is
defined to be
\begin{eqnarray*}
\Lambda^K_0(\psi) := \sum_{(H)\in \consub(K)} L^{\Z WH}(\psi^H,\psi^{>H})\cdot
[K/H] \in A(K)=U^K(*).
\end{eqnarray*}
We call the injective ring homomorphism
\[
\ch^K_0\colon A(K) \to  \bigoplus_{(H)\in \consub(K)} \Z, S  \mapsto  \{ |S^H| \}_{(H)\in \consub(K)}
\]
the \bet{character map}. Let $V$ be a finite-dimensional $K$-representation and let $\psi\colon V^c\to
V^c$ be a $K$-endomorphism of the one-point compactification
$V^c$. Define the \bet{equivariant degree} of $\psi$ to be
\[
\Deg^K_0(\psi):=\bigl(\Lambda^K_0(\psi)-1\bigr)\bigl(\Lambda^K_0(\id_{V^c})-1\bigr) \in A(K)=U^K(*).
\]
Let $G$ be a discrete group, $X$ a $G$-space and $f\colon X\to X$ an
equivariant endomorphism. Let $x\in X$ be a fixed point of $f$. We define
\begin{eqnarray*}
\Lambda_G(x,f) \colon \qquad \qquad U^G(G/G_x) & \to & \Lambda_G(X,f)\\
{}1\cdot [\tau \colon G/L\to G/G_x] & \mapsto & \overline{1_{x\circ \tau}} \cdot [x\circ \tau\colon G/L\to X],
\end{eqnarray*}
with $\overline{1_{x\circ \tau}}\in \Z \pi_1(X^L(x\circ \tau),x\circ
\tau)_{\phi_{x\circ \tau,\cst}'}$. Choosing $z\cong x\circ \tau\colon G/L\to X$ as a fixed
representative of the isomorphism class, this is
\[
\Lambda_G(x,f):  U^G(G/G_x) \to \Lambda_G(X,f), 1\cdot [\tau\colon G/L\to G/G_x]\mapsto \overline{\alpha_{x\circ \tau}} \cdot [z],
\]
where $\overline{\alpha_{x\circ \tau}}=(\sigma,[t])^*_{w,\cst}(\overline{1_{x\circ \tau}}) = \overline{ v f(t^{-1})
  \sigma^*(1_{x\circ \tau}) t} = \overline{ v f(t^{-1}) t }$ for an
isomorphism $(\sigma,[t])\in \Mor( z, x\circ \tau)$ and
$w=(\id,[v])\in \Mor(f(z),z)$.

\begin{Def}\label{4.6*}
Let $G$ be a discrete group, let $M$ be a
  cocompact smooth proper $G$-manifold and let $f\colon M\to M$ be a
  $G$-equivariant endomorphism such that $\Fix(f)\cap \partial M =
  \emptyset$ and such that for every $x\in \Fix(f)$ the determinant of the
  map $(\id_{T_x M} - T_x f)$ is different from zero. Then $G\setminus
  \Fix(f)$ is finite.
We define the
  \bet{generalized local equivariant Lefschetz
    class}\index{generalized local equivariant Lefschetz class} to be
\begin{eqnarray*}
\lambda^{loc}_G(f) & := & \sum_{Gx\in G\setminus \Fix(f)}
\Lambda_G(x,f) \circ \ind{G_x\subseteq
  G}\Bigl(\Deg^{G_x}_0\bigl((\id_{T_x M}-T_xf)^c\bigr) \Bigr).
\end{eqnarray*}
\end{Def}

\begin{Def}
We define the \bet{character map} $\ch_G(X,f)$ by
\begin{eqnarray*}
\Lambda_G(X,f) & \to & \bigoplus_{\overline{y}\in
  \Is\Pi(G,X)} \Q \pi_1(X^K(y),y)_{\phi_{y,w}'}\\
\bigl(\sum_{\overline{\alpha}} n_{\overline{\alpha}}\cdot
\overline{\alpha}\bigr)\cdot\overline{x} & \mapsto &
\sum_{\overline{y}} \sum_{\Aut(y)\cdot (\sigma, [t])\in \atop \Aut(y)\setminus \Mor(y,x)}  \hspace{-1em} \bigl|\bigl(\Aut(y)\bigr)_{(\sigma, [t])}\bigr|^{-1} (\sigma,[t])^* \Bigl(\sum_{\overline{\alpha}} n_{\overline{\alpha}}\cdot \overline{\alpha}\Bigr).
\end{eqnarray*}
\end{Def}

This character map is a generalization of
\cite[Definition~5.1]{lueck-rosenberg}.

\begin{Lem}\label{5.3*}
The character map $\ch_G(X,f)$ is injective.
\end{Lem}
\begin{proof}
Let $u=\sum_{i=1}^{k} a_i \cdot \overline{x_i}\in \Lambda_G(X,f)$,
where $a_i\in \Lambda_G(X,f)_{{x_i}}$, with
  $\ch_G(X,f)(u)=0$. Let the $x_i$ be ordered in accordance with the
  partial ordering on $\Pi(G,X)$ given by $\overline{x}\leq
  \overline{y}\Leftrightarrow \Mor(x,y)\neq 0$, so $x_i\leq x_j
  \Rightarrow i\leq j$. Suppose without loss of generality that
  $a_k \neq 0$. If $\ch_G(X,f)(\overline{x_i})_{\overline{x_k}}\neq
    0$, then $\Mor(x_k,x_i)\neq 0$, which implies $x_k\leq x_i$ and
    thus $x_k=x_i$. Since $\ch_G(X,f)(a_k\cdot
    \overline{x_k})_{\overline{x_k}}= a_k$, we obtain
    $0=\ch^G(X,f)(u)=a_k$, a contradiction.
\end{proof}

We now
calculate the value of $\ch_G(X,f)$ on the generalized equivariant
Lefschetz invariant $\lambda_G(f)$ and on the generalized local
Lefschetz class $\lambda_G^{loc}(f)$, in analogy to~\cite[Lemma~5.4
and Lemma~5.9]{lueck-rosenberg}. Using Theorem~\ref{2.1*}, these values
will turn out to be equal, proving Theorem~\ref{0.2*}.

\begin{Lem}\label{5.4*}
Let $f\colon X\to X$ be a $G$-equivariant endomorphism of a finite proper $G$-CW-complex $X$. Let $\overline{y}$ be an isomorphism class of objects $y\colon G/K\to X$ in $\Pi(G,X)$. Then
\[
\ch_G(X,f)\bigl(\lambda_G(f)\bigr)_{\overline{y}}= L^{\Q \Aut(y)}\bigl(\widetilde{f^K(y)}\bigr).
\]
\end{Lem}

\begin{proof}
We first consider the case $X^K(f(y))= X^K(y)$. We write the
$p$-skeleton $X_p$ as a pushout and call
$x_{p,i}\colon G/H_i\to X$ for $0\leq
i\leq n_p$ the centers of the equivariant $p$-cells.

The $G$-CW-structure on $X$ induces an $\Aut(y)$-CW-structure on $\widetilde{X^K(y)}$. We obtain a pushout diagram of $\Aut(y)$-spaces
\[
\xymatrix@R=1pc{
\protect \coprod_{i=1}^{n_p} \Mor(y,x_{p,i})\times S^{p-1} \ar[r] \ar[d] & \protect\widetilde{X^K(y)_{p-1}} \ar[d]\\
\protect \coprod_{i=1}^{n_p} \Mor(y,x_{p,i})\times D^{p} \ar[r] & \protect\widetilde{X^K(y)_{p}}
}
\]
(If $p-1\leq 1$, we use the cover corresponding
to $\pi_1(X^K(y),y)$.) We denote the $p$-cells of $\widetilde{X^K(y)}$
by $e_{(\sigma,[t]),p,i}:=(\sigma,[t])\times
{\oD}{}^{p}$, where $(\sigma,[t])\in \Mor( y, x_{p,i})$. The $\Aut(y)$-orbit of the cell $e_{(\sigma,[t]),p,i}$ corresponds to the $\Aut(y)$-orbit of $(\sigma,[t])$, so we conclude from Lemma~\ref{1.9*} that
\begin{eqnarray*}
\lefteqn{L^{\Q \Aut(y)}\bigl(\widetilde{f^K(y)}\bigr)}\\
& = & \sum_{p\geq 0} (-1)^p \sum_{i=1}^{n_p}\hspace{-0.5em} \sum_{\Aut(y)\cdot
  (\sigma,[t]) \in \atop \Aut(y)\setminus \Mor(y,x_{p,i})}\hspace{-0.5em} \bigl|\Aut(y)_{(\sigma,[t])}\bigr|^{-1} \cdot \inc_{\phi_{y,w}}\bigl(\widetilde{f^K(y)}, e_{(\sigma,[t]),p,i}\bigr).
\end{eqnarray*}
Analogously, we have for any $x\colon G/H\to X$ a pushout diagram
\[
\xymatrix@R=1pc{
\protect \coprod_{i=1}^{n_p} \Mor(x,x_{p,i})\times S^{p-1} \ar[r] \ar[d] & \protect\widetilde{X^H(x)_{p-1}\cup X^{>H}(x)} \ar[d] \\
\protect \coprod_{i=1}^{n_p} \Mor(x,x_{p,i})\times D^{p} \ar[r] & \protect\widetilde{X^H(x)_{p}\cup X^{>H}(x)} .
}
\]
Lemma~\ref{1.9*} yields
\[
L^{\Z \Aut(x)}\bigl(\widetilde{f^H(x)},\widetilde{f^{>H}(x)}\bigr)
 =  \sum_{p\geq 0} (-1)^p \hspace{-0.5em} \sum_{i=1 \atop \overline{x_{p,i}}=\overline{x}}^{n_p} \hspace{-0.5em} (\tau,[s])^* \inc_{\phi_{x_{p,i},\cst}}\bigl(\widetilde{f^H(x_{p,i})},e_{p,i}\bigr),
\]
for $(\tau,[s])\in \Mor(x,x_{p,i})$ any morphism. Inserting this formula into the definition of $\ch_G(X,f) \bigl(\lambda_G(f)\bigr)_{\overline{y}}$ proves the claim.

Now we consider the case that $X^K(f(y))\neq X^K(y)$. This implies that
$X^H(f(x))\neq X^H(x)$ for all $\overline{x}$ with $\Mor(y,x)\neq \emptyset$, so
$\lambda_G(f)_{\overline{x}}=0$ for all $\overline{x}$ with
$\Mor(y,x)\neq \emptyset$. Therefore $\ch_G(X,f)(\lambda_G(f))_{\overline{y}}=0$.
\end{proof}

\begin{Lem} \label{5.9*}
Let $G$ be a discrete group and let $M$ be a cocompact smooth proper $G$-manifold. Let $f\colon M\to
M$ be a smooth $G$-equivariant map. Suppose that $\Fix(f)\cap \partial M =
  \emptyset$ and that for any $x\in
\Fix(f)$ the determinant $\det(\id_{T_x M}-T_x f)$ is different from
zero. Then the set $G\setminus \Fix(f)$ is finite. Let $y\colon G/K\to M$ be an object in $\Pi(G,M)$. Then the set
$WK_y\setminus \Fix\bigl(f|_{M^K(y)}\bigr)$ is finite and we get
\begin{eqnarray*}
\lefteqn{\ch_G(M,f)\bigl(\lambda^{loc}_G(f)\bigr)_{\overline{y}}}\\
& = & \sum_{WK_y \cdot x\in \atop WK_y \setminus \Fix (f|_{M^K(y)} )}
|(WK_y)_x|^{-1} \deg \Bigl( \bigr(\id_{T_x M^K(y)} - T_x(f|_{M^K(y)})\bigr)^c  \Bigr)\cdot \overline{\alpha_x},
\end{eqnarray*}
where $\overline{\alpha_x}=\overline{v f(t^{-1}) t}\in
\pi_1(X^K(y),y)_{\phi'}$ for $(\sigma,[t])\in \Mor(y,x)$ and $w=(\id,[v])\in \Mor(f(y),y)$.
\end{Lem}

\begin{proof}
The set $G\setminus \Fix(f)$ is finite since $M$ is cocompact and the
fixed points are isolated. Analogously, $G\setminus GM^K(y) = WK_y
\setminus M^K(y)$ is compact with isolated fixed points, so
$WK_y\setminus \Fix(f|_{M^K(y)})$ is finite. Let $x\colon G/G_x\to M$
be a fixed point of $f$.

We first show that for each $u\in U^{G_x}(*)$
we have
\[
{\ch_G(M,f)}_{\overline{y}}\Lambda_G(x,f) \ind{G_x\subseteq G}(u)
\hspace{-0.3em} = \hspace{-2.5em} \sum_{\Aut(y)\cdot (\sigma,[t])\in \atop \Aut(y)\setminus \Mor(y,x)}\hspace{-1.8em}  \bigl|\Aut(y)_{(\sigma,[t])}\bigr|^{-1}  (\sigma,[t])^* \overline{1_{x}} \ch^{G_x}_{0}(u)_{(K_\sigma)}
\]
where $\overline{1_{x}}\in \Z\pi_1(X^{G_x}(x),x)_{\phi_{x,\cst}}$. Here $(K_\sigma)=(g_{\sigma}^{-1}K
g_{\sigma}) \in \consub(G_x)$ for $\sigma:G/K\to G/G_x, gK\mapsto g
g_\sigma G_x$. Let $u=[G_x/L]\in U^{G_x}(*)$ be a basis element and $\pr\colon G/L\to G/G_x$ the projection. Then
\[
{\ch_G(M,f)}_{\overline{y}} \Lambda_G(x,f) \ind{G_x\subseteq G}\bigl([G_x/L]\bigr) = \hspace{-3em} \sum_{\Aut(y)\cdot (\tau,[t])\in \atop \Aut(y)\setminus \Mor(y,x\circ \pr)} \hspace{-2em} \bigl|\Aut(y)_{(\tau,[t])}\bigr|^{-1} (\tau,[t])^*\overline{1_{x\circ \pr}}
\]
where $\overline{1_{x\circ \pr}}\in
\Z\pi_1(X^L(x\circ \pr),x\circ \pr)_{\phi_{x\circ \pr,\cst}}$.

Defining $q\colon \Mor(y,x\circ \pr)\to
\Mor(y,x), (\tau,[t])\mapsto (\pr\circ\tau,[t])$ we have
\[
\Mor(y,x\circ \pr) = \hspace{-0.5em}\coprod_{(\sigma,[t])\in \atop \Mor(y,x)}
\hspace{-0.5em}q^{-1}(\sigma,[t]) = \hspace{-0.5em}\hspace{-0.5em}\coprod_{\Aut(y)\cdot (\sigma,[t])\in \atop \Aut(y)\setminus \Mor(y,x)} \hspace{-1em}\Aut(y)\times_{\Aut(y)_{(\sigma,[t])}}q^{-1}(\sigma,[t]).
\]
The $\Aut(y)_{(\sigma,[t])}$-set $q^{-1}(\sigma,[t])=\coprod_{i\in I(\sigma,[t])} \Aut(y)_{(\sigma,[t])}/A_i$ is a finite disjoint union of orbits, thus we have a bijection of $\Aut(y)$-sets
\[
\Mor(y,x\circ \pr)=\coprod_{\Aut(y)\cdot (\sigma,[t])\in \atop
  \Aut(y)\setminus \Mor(y,x)} \quad \coprod_{i\in I(\sigma,[t])} \Aut(y)/A_i.
\]
We know for $(\tau,[t])\in \Mor(y,x\circ \pr)$ that $(\tau,[t])^* \overline{1_{x\circ \pr}} = \overline{v
  f(t^{-1}) t } = (\pr \circ \tau, [t])^* \overline{1_x}$. An orbit
$\Aut(y)/A_i$ corresponds to exactly one orbit $\Aut(y)\cdot
(\tau,[t])$, where $i\in I(\pr \circ \tau,[t])$, so $|\Aut(y)_{(\tau,[t])}|=|A_i|$, hence
\[
\sum_{\Aut(y)\cdot (\tau,[t])\in \atop \Aut(y)\setminus
    \Mor(y,x\circ\pr)}\hspace{-1pc} \hspace{-1.4em} \bigl|\Aut(y)_{(\tau, [t])}\bigr|^{-1} (\tau,[t])^* \overline{1_{x\circ \pr}}
= \hspace{-1pc} \hspace{-0.5em}\sum_{\Aut(y)\cdot (\sigma,[t])\in \atop \Aut(y)\setminus \Mor(y,x)} \sum_{i\in I(\sigma,[t])}  \hspace{-1em} |A_i|^{-1} (\sigma,[t])^* \overline{1_{x}}.
\]
We have
$
\bigl|q^{-1}(\sigma,[t])\bigr|=\bigl|\Aut(y)_{(\sigma,[t])}\bigr|\cdot \sum_{i\in I(\sigma,[t])}|A_i|^{-1}.
$
Since $q$ does not change the $[t]$-part, as in~\cite[Lemma~5.9, Equation~5.14]{lueck-rosenberg} we have
$
\bigl|q^{-1}(\sigma,[t])\bigr| = \bigl|G_x/L^{g_\sigma^{-1}Kg_{\sigma}}\bigr| = \ch_0^{G_x}\bigl([G_x/L]\bigr)_{(K_\sigma)}.
$

Inserting these equations into the above formula, we obtain the desired equation for all $[G_x/L]\in
U^{G_x}(*)$, thus for all $u\in U^{G_x}(*)$.

We know~\cite[Equations~5.16 and 5.17]{lueck-rosenberg}
\[
\ch_0^{G_x}\bigl(\Deg^{G_x}_0\bigl((\id_{T_x M} - T_x
f)^c\bigr)\bigr)_{(K_\sigma)} \hspace{-1em} = \deg\bigl(\bigl(\id_{T_{g_{\sigma}x} M^K(y)}- T_{g_{\sigma}x}(f|_{M^K(y)})\bigr)^c\bigr).
\]
We have
$
\coprod_{G\setminus\Fix(f)} \Aut(y)\setminus \Mor(y,x) \cong WK_y\setminus \Fix\bigl(f|_{M^K(y)}\bigr).
$
Under this bijection, $\Aut(y)\cdot (\sigma,[t])\mapsto WK_y \cdot
t(0)$, where $t(0)=x\circ \sigma(1K)=x(g_{\sigma} G_x)=g_{\sigma}x(1
G_x)=g_{\sigma}x$ for $\sigma\colon G/K\to G/G_x, g K\mapsto g g_{\sigma} G_x$. Since $|\Aut(y)_{(\sigma,[t])}|=|(WK_y)_{g_{\sigma}x}|$, inserting the above results into the formula for $\ch_G(M,f)\bigl(\lambda^{loc}_G(f)\bigr)_{\overline{y}}$ yields the claim.
\end{proof}

The final ingredient in the proof of the Theorem~\ref{0.2*} is a refinement of the orbifold Lefschetz fixed point
theorem~\cite[Theorem~2.1]{lueck-rosenberg}.

\begin{Thm}\label{2.1*}
Let $G$ be a discrete group extension $1\to \pi \to G \to W \to 1$ with
endomorphism $\phi\colon G\to G$ such that
$\phi_W=\id_{W}$. Let $M$ be a connected simply connected
cocompact proper $G$-manifold such that $\pi$ operates freely on $M$,
and let $f\colon M\to M$ be a smooth $\phi$-twisted map. Denote by
$\overline{f}\colon \overline{M}\to \overline{M}$ the $W$-equivariant
map induced on the manifold $\overline{M}:=\pi\setminus M$ by dividing out the
$\pi$-action. Suppose that $\Fix(\overline{f})\cap \partial
\overline{M}=\emptyset$ and that for every $x\in \Fix(\overline{f})$
the determinant of the map~$(\id_{T_x\overline{M}} - T_x\overline{f})$
is different from zero. Then $W\setminus \Fix(\overline{f})$ is
finite, and
\[
L^{\Q G}(f)=\sum_{W\cdot x \in W\setminus \Fix(\overline{f})}
|W_x|^{-1} \cdot \deg \bigl((\id_{T_x\overline{M}} -
T_x\overline{f})^c\bigr) \cdot \overline{\alpha_x}.
\]
Here $(\id_{T_x\overline{M}} -
T_x\overline{f})^c\colon (T_x\overline{M})^c \to (T_x\overline{M})^c$
is an endomorphism of the one-point compactification of $T_x\overline{M}$.
\end{Thm}

\begin{proof}
The proof is analogous to the proof of the orbifold Lefschetz fixed
point theorem~\cite[Theorem~2.1]{lueck-rosenberg}. We apply the
construction from that proof to the $W$-equivariant map $\overline{f}\colon
\overline{M}\to \overline{M}$. We obtain a $W$-equivariant map
$\overline{f'}$ such that $\Fix(\overline{f})=\Fix(\overline{f'})$,
the map $\overline{f'}$ is of the desired form around the fixed
points and agrees with $\overline{f}$ outside a neighborhood of
the fixed points. The desired form is that
$\exp_{x,\varepsilon_1}^{-1} \circ \overline{f} \circ \exp_{x,\varepsilon_2}$
and $\overline{f'}$ agree on $D_{\varepsilon}T_x \overline{M}$ for
some $\varepsilon>0$ and for all
$x\in \Fix(\overline{f})$. Here $\exp_{x,\varepsilon}\colon
D_{\varepsilon} T_x \overline{M} \xrightarrow{\sim} N_{x,\varepsilon}$ denotes the
exponential map.

We lift this construction to $M$ and the
$\phi$-twisted endomorphism $f\colon M\to M$ by extending
$\phi$-twisted $G$-equivariantly: Let $z$ be a lift of a fixed
point~$x$, then $f(z)=\alpha_z \cdot z$ with $\alpha_z \in \pi$, and
we set
$f'|_{U_z}:=\alpha_z \cdot \varphi \overline{f'} \varphi^{-1}$ on a
neighborhood of $z$ that is isomorphic to a neighborhood of $x$ via
the isomorphism $\varphi\colon U_z\xrightarrow{\sim} U_x$ coming from
the covering map. Around another lift $\beta\cdot z$, we set
$f'|_{U_{\beta z}}:=\phi(\beta)\cdot f'|_{U_z} \cdot \beta^{-1}$ on a
neighborhood of $\beta z$. (Note that $\alpha_{\beta
  z}=\phi(\beta)\cdot \alpha_z \cdot \beta^{-1}$.) We obtain a
$\phi$-twisted map $f'\colon M \to M$ such that $\Fix(f)=\Fix(f')$,
the map $f'$ is of the desired form around the
orbits of the fixed points, and it agrees with $f$ outside a
neighborhood of the orbits of the fixed points.

Analogously, we lift further constructions such as the
$W$-equivariant triangulations $K'(\overline{M})$,
$K''(\overline{M})$ and the $W$-homotopy $\overline{h}$ from $\overline{f}$ to $\overline{h_1}$. The construction of
$K''(\overline{M})$ can be done such that there is at most one
fixed point $x$ of $\overline{h_1}$ in each $\overline{e}\in
K''(\overline{M})$. We denote $\overline{h_1}$ again by
$\overline{f}$. Then
\[
\inc_\phi(f,e)=\begin{cases}
                   \inc(\overline{f},\overline{e}) \cdot \overline{\alpha_x} & \text{ if there is a fixed point } x\in \overline{e}\\
                   0 & \text{ else}.
               \end{cases}
\]
Here for a basepoint $y$ of $\overline{M}$, a path $v$ from $y$ to
$\overline{f}(y)$ and a morphism $(\sigma,[t])\colon y\to x$ from $y$
to $x$ we set $\overline{\alpha_x}:=\overline{v f(t^{-1}) t}$. The direction
of $v$ corresponds to our
usual convention that $w=(\id,[v])\in \Mor(\overline{f}(y),y)$. Note that
$\overline{\alpha_x}$ does not depend on the choice of $(\sigma,[t])$
because of Lemma~\ref{lemindepofsigmat}. We see that
$|G_e|=|W_{\overline{e}}|$ because $\pi$ operates freely. So using the
construction of~\cite[proof of Theorem~2.1]{lueck-rosenberg} applied to the $W$-equivariant
map $\overline{f}$, where $\{x_1,\ldots,x_k\}$ is a complete set of
representatives of $W$-orbits of fixed points of $\overline{f}$, we
obtain
\begin{eqnarray*}
L^{\Q G}(f) = L^{\Q G}(h_1)
& = & \sum_{p\geq 0}(-1)^p \sum_{G\cdot e\in G\setminus I_p(K''(M))} |G_e|^{-1} \inc_{\phi}(f,e) \\
& = & \sum_{p\geq 0}(-1)^p \sum_{W\cdot \overline{e}\in W\setminus I_p(K''(\overline{M}))} |W_{\overline{e}}|^{-1} \inc(\overline{f},\overline{e})\cdot \overline{\alpha_{x_{\overline{e}}}}\\
& = & \sum_{i=1}^{k} |W_{x_i}|^{-1} \deg \bigl((\id - T_{x_i} \overline{f})^c\bigr) \cdot \overline{\alpha_{x_i}}\\
& = & \sum_{W\cdot x \in W \setminus \Fix(\overline{f})} |W_x|^{-1} \cdot \deg\bigl((\id - T_{x}\overline{f})^c\bigr) \cdot \overline{\alpha_x}.
\end{eqnarray*}
\vskip-6ex
\end{proof}

We have assembled all information necessary for the proof of the
refined equivariant Lefschetz fixed point theorem.

\begin{Thmb}
Let $G$ be a discrete group, let $M$ be a
  cocompact proper smooth $G$-manifold and let $f\colon M\to M$ be a
  $G$-equivariant endomorphism such that $\Fix(f)\cap \partial M =
  \emptyset$ and such that for every $x\in \Fix(f)$ the determinant of the
  map $(\id_{T_x M} - T_x f)$ is different from zero. Then
\[
\lambda_G(f)=\lambda^{loc}_G(f).
\]
\end{Thmb}

\begin{proof}
Using Lemma~\ref{5.4*}, Theorem~\ref{2.1*} and Lemma~\ref{5.9*} we have
\begin{eqnarray*}
\lefteqn{\ch_G(M,f)\bigl(\lambda_G(f)\bigr)_{\overline{y}}}\\
& = & L^{\Q \Aut(y)}\bigl(\widetilde{f^K(y)}\bigr)\\
& = & \sum_{WK_y \cdot x\in \atop WK_y \setminus
  \Fix (f|_{M^K(y)} )} \bigl|(WK_y)_x\bigr|^{-1} \deg\Bigl(
\bigl(\id_{T_x M^K(y)} - T_x(f|_{M^K(y)})\bigr)^c  \Bigr) \cdot \overline{\alpha_x}\\
& = & \ch_G(M,f)\bigl(\lambda^{loc}_G(f)\bigr)_{\overline{y}}
\end{eqnarray*}
for all $\overline{y}\in \Is \Pi(G,X)$. By injectivity of $\ch_G(M,f)$,
shown in Lemma~\ref{5.3*}, we
obtain the refined equivariant Lefschetz fixed point theorem.
\end{proof}

From the fact that the units of the Burnside ring $A(K)$ are only
$\{1,-1\}$ if $K$ is a finite group of odd order~\cite[Proposition~1.5.1]{tomdieck79}, we obtain in analogy
to~\cite[Example~4.7]{lueck-rosenberg} the following example.

\begin{Ex}
Let $G$ be a discrete group and let $M$ be a cocompact proper smooth
$G$-manifold. Suppose that all isotropy groups $G_x$ of points $x\in
M$ are of odd order. Let $f\colon M\to M$ be a smooth $G$-equivariant
map such that $\Fix(f)\cap \partial M = \emptyset$ and such that for
every $x\in \Fix(f)$ the determinant of the map $(\id_{T_x M} - T_x
f)$ is different from zero. Then
\[
\lambda^{loc}_G(f) = \sum_{Gx\in G\setminus \Fix(f)} \frac{\det\bigl(1-T_x(f)\colon T_x(M)\to T_x(M)\bigr)}{\left|\det\bigl(1-T_x(f)\colon T_x(M)\to T_x(M)\bigr)\right|} \cdot \overline{\alpha_x}  ,
\]
where for an isomorphism $(\sigma,[t])\in
\Mor(z,x)$ and a path $w=(\id,[v])\in \Mor(f(z),z)$ we have $\alpha_x=\overline{v f(t^{-1})t}\in \Z
\pi_1(X^{G_x}(z),z)_{\phi'}$.
\end{Ex}



\def\cprime{$'$}

\end{document}